\documentclass{amsart}

\usepackage{amscd,amssymb,amsmath}
\usepackage{mathrsfs}
\usepackage[all]{xy}

\newcommand\mylabel[1]{\label{#1}}

\newtheorem{theorem}{Theorem}[section]
\newtheorem*{maintheorem}{Theorem}
\newtheorem{lemma}[theorem]{Lemma}
\newtheorem{proposition}[theorem]{Proposition}
\newtheorem{corollary}[theorem]{Corollary}

\theoremstyle{definition}
\newtheorem{definition}[theorem]{Definition}

\newtheorem*{acknowledgement}{Acknowledgement}
\newtheorem*{remarkintro}{Remark}

\theoremstyle{remark}

\DeclareFontFamily{U}{wncy}{}
\DeclareFontShape{U}{wncy}{m}{n}{<->wncyr10}{}
\DeclareSymbolFont{mcy}{U}{wncy}{m}{n}
\DeclareMathSymbol{\Sh}{\mathord}{mcy}{"58}

\newcommand{\ZZ}	{\mathbb{Z}}
\newcommand{\QQ}	{\mathbb{Q}}

\newcommand{\idealq}    {\mathfrak{q}}
\newcommand{\ideala}    {\mathfrak{a}}
\newcommand{\idealb}    {\mathfrak{b}}

\newcommand  {\shE}     {\mathscr{E}}

\newcommand  {\shT}     {\mathscr{T}}

\newcommand  {\shU}     {\mathscr{U}}

\newcommand  {\calC}     {\mathcal{C}}
\newcommand  {\calD}     {\mathcal{D}}

\newcommand  {\calF}     {\mathcal{F}}

\newcommand  {\calS}     {\mathcal{S}}

\newcommand  {\bg}       {{\text{\rm big}}}

\newcommand  {\can}     {{\rm \text{can}}}
\newcommand  {\Card}    {\operatorname{Card}}

\newcommand  {\Cat}     {{\text{\rm Cat}}}

\newcommand  {\Cov}     {\operatorname{Cov}}

\newcommand  {\et}      {{\text{\rm et}}}

\newcommand  {\fppf}    {{\text{{\rm fppf}}}}

\newcommand  {\fpuo}    {{\text{{\rm fpuo}}}}

\newcommand  {\Grp}     {\operatorname{Grp}}

\newcommand  {\Hom}     {\operatorname{Hom}}

\newcommand  {\id}      {{\operatorname{id}}}

\newcommand  {\Loc}     {\operatorname{Loc}}
\newcommand  {\dirlim}  {\varinjlim}
\newcommand  {\invlim}  {\varprojlim}

\newcommand  {\lra}     {\longrightarrow}

\newcommand  {\maxid}   {\mathfrak{m}}

\newcommand  {\Nil}     {\operatorname{Nil}}

\newcommand  {\primid}  {\mathfrak{p}}
\renewcommand{\O}       {\mathscr{O}}

\newcommand  {\op}      {{\operatorname{op}}}

\newcommand  {\Points}  {\operatorname{Points}}

\newcommand  {\pr}      {\operatorname{pr}}

\newcommand  {\PSh}     {\operatorname{PSh}}

\newcommand  {\quadand} {\quad\text{and}\quad}

\newcommand  {\ra}      {\rightarrow}

\newcommand  {\res}     {\operatorname{res}}

\newcommand  {\Set}     {{\text{\rm Set}}}
\newcommand  {\Sch}     {{\text{\rm Sch}}}

\renewcommand{\Sh}      {{\operatorname{Sh}}} 

\newcommand  {\sob}     {\text{{\rm sob}}}

\newcommand  {\Spec}    {\operatorname{Spec}}

\newcommand {\Zar}      {{\text{\rm Zar}}}

\def\mydate{\number\day\space\ifcase\month \or January\or February\or March\or 
April\or May\or June\or July\or
August\or September\or October\or November\or December\fi \space\number\year}

\DeclareFontFamily{U}{wncy}{}
\DeclareFontShape{U}{wncy}{m}{n}{<->wncyr10}{}
\DeclareSymbolFont{mcy}{U}{wncy}{m}{n}
\DeclareMathSymbol{\SH}{\mathord}{mcy}{"58}

\begin{document}

\title[Points in the fppf topology]
      {Points in the fppf topology}

\author[Stefan Schr\"oer]{Stefan Schr\"oer}
\address{Mathematisches Institut, Heinrich-Heine-Universit\"at,
40225 D\"usseldorf, Germany}
\curraddr{}
\email{schroeer@math.uni-duesseldorf.de}


\dedicatory{Final version, 22 January 2016}

\begin{abstract}
Using methods from commutative algebra and topos-theory, we construct topos-theoretical points  for the fppf topology of a scheme. These points are
indexed by both a geometric point and a limit ordinal. 
The resulting stalks of the structure sheaf are what we call fppf-local rings.
We show that for such rings  all   localizations at primes
are henselian with algebraically closed residue field, and relate
them to AIC and TIC rings.
Furthermore, we give an abstract criterion ensuring that  two sites have point spaces
 with identical sobrification. This applies in particular
to some standard Grothendieck topologies considered in algebraic geometry:
Zariski, \'etale, syntomic, and fppf.
\end{abstract}

\maketitle
\tableofcontents

\section*{Introduction}
One of the major steps in Grothendieck's program  to prove the Weil Conjectures 
was the introduction of topoi \cite{SGA 4a}, thus lying the foundations for   \'etale cohomology. Roughly speaking, a \emph{topos} is
a category $\shE$ that is equivalent to the category of sheaves $\Sh(\calC)$ on some \emph{site} $\calC$.
The latter is a category endowed with a Grothendieck topology, which gives
the objects   a role similar to the  open subsets of a topological space.
(Set-theoretical issues will be neglected in the introduction, but treated with care in what follows.)

One may perhaps say that topoi are the true incarnation of our  notion of space,
keeping exactly what is necessary to pass back and forth between local and global,
to apply geometric intuition, and to use   cohomology. 
This comes, of course, at the price  of erecting a frightening technical apparatus.
Half-way between  sites and topological spaces dwell the so-called \emph{locales},
sometimes referred to as \emph{pointless spaces}, which are certain ordered sets $L$ having analogous order properties
like the collection  $\shT$ comprising the open subsets from a topological space  $(X,\shT)$.

A considerable part of \cite{SGA 4a} deals with the notion of points for a topos.
Roughly speaking,  \emph{topos-theoretical points}  are continuous maps of topoi $P:(\Set)\ra\shE$,
where the category of sets is regarded as the topos of sheaves on a singleton space.
Such a continuous map consists of a \emph{stalk functor} $P^{-1}:\shE\ra(\Set)$, which must commute with finite inverse limits, and 
a direct image functor $P_*:(\Set)\ra\shE$, related by an adjunction.
After one chooses an equivalence  $\shE\simeq\Sh(\calC)$,   
the points  $P$ can be recovered via certain  pro-objects $(U_i)_{i\in I}$ of \emph{neighborhoods}
$U_i\in\calC$, such that $F_P=P^{-1}(F)=\dirlim_i\Gamma(U_i,F)$ for all sheaves $F$ on the site $\calC$.

The topos-theoretical points $P$ form the \emph{category of points} $\Points(\shE)$, and their isomorphism classes $[P]$
comprise the \emph{space of points}   $|\shE|$, where the topology comes from the subobjects  of the terminal object   $e\in\shE$.
In the special case that $\shE=X_\Zar=\Sh(X)$ is the topos of sheaves on a topological space
$X$,   there is  continuous map $X\ra|\shE|$,
which  can be identified with the \emph{sobrification} $X\ra X_\sob$, in other words,
the universal map into a topological space where each irreducible closed subset has a unique generic points.
In particular, for schemes $X$ endowed with the Zariski topology we have an identification $X=|\shE|$.

Another important result is \emph{Deligne's Theorem},   which asserts that topoi fulfilling certain technical
finiteness conditions have \emph{enough points}, that is,
stalk functors detect monomorphisms (\cite{SGA 4b}, Appendix to Expose VI). 
Note, however,  that there are   examples of topoi having no point at all,
and examples of topoi having   ``large'' spaces of points (\cite{SGA 4a}, Expose IV, Section 7).

The goal of this paper is to investigate the space of points $|\shE|$ for various topoi
occurring in algebraic geometry, in particular for the \emph{fppf topology}.
The fppf topos was studied, for example,  by Milne (\cite{Milne 1980}, Chapter III, \S3 and \cite{Milne 1986}, Chapter III),
Shatz (\cite{Shatz 1972}, Chapter VI and \cite{Shatz 1964}, \cite{Shatz 1966}, \cite{Shatz 1968}),
Waterhouse \cite{Waterhouse 1975}, and the Stacks Project \cite{Stacks Project}.
Somewhat surprisingly,  very little seems to be known about the points.
Interesting result on several other Grothendieck topologies
were recently obtained by Gabber and Kelly \cite{Gabber; Kelly 2015}.

Note that one ususally defines the fppf topos via sheaves on the ``big'' site $(\Sch/X)$
of all $X$-schemes. In this paper, however, we will mainly consider the topos obtained
by sheaves on the ``small'' site $(\fppf/X)$ comprising flat schemes that are locally of finite presentation.

To gain flexibility and facilitate applications,
it seems preferable to work in an axiomatized situation, having nevertheless the fppf topos in mind.
One of our main result is a sufficent criterion for   adjoint functors between two abstract sites 
\begin{equation}
\label{adjoint pair}
u:\calC_f\ra\calC_z\quadand v:\calC_z\ra\calC_f
\end{equation}
to induce  homeomorphisms between locales:

\begin{maintheorem}[see Thm.\ \ref{locale isomorphic}]
Suppose the adjoint functors  above satisfy the conditions (TL 1) -- (TL 4) given in Section \ref{topoi same locales}.
Then the induced continuous maps of locales $\epsilon:\Loc(\shE_f)\ra\Loc(\shE_z)$
is a homeomorphism. Moreover,  we get an embedding 
$|\shE_f|_\sob\subset |\shE_z|_\sob$ of sober spaces. The latter is an equality if
the map $|\shE_f|\ra|\shE_z|$ admits a section.
\end{maintheorem}

Note that these conditions are rather technical to formulate, but easy to verify in practice.
The result applies  for the fppf site $\calC_f=(\fppf/X)$ comprising
flat $X$-schemes that are locally of finite presentation, 
and the Zariski site $\calC_z=(\Zar/X)$ of open subschemes,
and their corresponding topoi $\shE_f=X_\fppf$ and $\shE_z=X_\Zar$. 
The morphisms in the above sites are the $X$-morphisms and inclusion maps, respectively.
Since we are able to construct   sections for the canonical map $|X_\fppf|\ra|X_\Zar|$, we get:

\begin{maintheorem}[see Thm.\ \ref{identifiaction sober}]
For every scheme $X$,  the continuous map of topoi $X_\fppf\ra X_\Zar$
induces an identification
$|X_\fppf|_\sob=|X_\Zar|_\sob=X$ of sober spaces.
\end{maintheorem}

Similar results hold for the Nisnevich topology, the \'etale topology, and the syntomic topology.
Note that for the \'etale site, the much stronger result $|X_\et|=X$ is true by \cite{Stacks Project}, Lemma 44.29.12. 
This, however, becomes false in the fppf topology, as we shall see below. 
Unfortunately, our method, as it stands, does not   apply to the fpqc site, because 
fpqc morphisms are not necessarily open maps.

In order to construct explicit points $P:(\Set)\ra(X_\fppf)$ for the fppf topos, 
we introduce the notion of \emph{fppf-local rings}, which are local rings $R$ for which
any fppf algebra $A$ admits a retraction, in other words, the morphisms $\Spec(A)\ra\Spec(R)$ has a section.
Such rings should be regarded as generalizations of algebraically closed fields. 
However, they have highly unusual properties from the point of view of commuative algebra. For example, their 
formal completion $\hat{R}=\invlim_n R/\maxid^n $ coincides with the residue field $\kappa=R/\maxid$.\
Rings with similar properties were studied by Gabber and Romero \cite{Gabber; Ramero 2003}, in the contex of
``almost mathematics''. Here is another amazing property:

\begin{maintheorem}[see Thm.\ \ref{amazing property}]
If $R$ is fppf-local, then the local rings $R_\primid$ are strictly local with
algebraically closed residue field, for all prime ideals $\primid\subset R$.
\end{maintheorem}

Moreover, we   relate fppf-local rings to the so-called
\emph{TIC rings} introduced by  Enochs \cite{Enochs 1968} and further studied by Hochster \cite{Hochster 1970}, 
and the \emph{AIC rings} considered by Artin \cite{Artin 1971}. Furthermore, we
show that the  stalks $\O_{X_\fppf,P}$ of the   structure sheaf at topos-theoretical points are examples
of fppf-local rings. Throughout, the term \emph{strictly local} denotes local henselian rings with
separably closed residue fields.

We then use ideas of Picavet \cite{Picavet 2007} to  construct, for each strictly local ring $R$
and each limit ordinal $\lambda$, some faithfully flat, integral ring extension $R_\lambda$ that is fppf-local.
Roughly speaking, the idea  is to form the tensor product over ``all'' finite fppf algebras,
and to iterate  this via tranfinite recursion, until reching the limit ordinal $\lambda$.
Note that there is a close analogy to the Steinitz's original construction of   algebraic closures
for fields (\cite{Steinitz 1910}, Chapter III).

This is next used to produce, for each geometric point $\bar{a}:\Spec(\Omega)\ra X$ on a scheme $X$ and each limit
ordinal $\lambda$, a topos-theoretical point 
$$
P=P_{\bar{a},\lambda}:(\Set)\ra X_\fppf
$$ 
with $\O_{X_\fppf,P}=(\O_{X,\bar{a}})_\lambda$. Here the main step is to construct a suitable \emph{pro-object}
  $(U_i)_{i\in I}$ of flat $X$-schemes locally of finite presentation yielding the stalk functor.
The index category $I$ will consists of certain 5-tuples of $X$-schemes and morphisms between them.
This   gives the desired section:

\begin{maintheorem}[see Thm.\ \ref{continuous section}]
For each limit ordinal $\lambda$, the  map $a\mapsto P_{\bar{a},\lambda}$
induces a continuous section for the canonical map  $|X_\fppf|\ra X$.
\end{maintheorem}

The paper is organized as follows:
In Section \ref{universes sites topoi}, we recall some basic definitions and result from topos-theory,
also paying special attention to set-theoretical issuses. 
Following Grothendieck, we  avoid the use of the ambiguous notion of ``classes'', and use universes instead.
Section \ref{topoi same locales} contains the sufficient criterion that  
two sites have homeomorphic locales and sobrified spaces.
This is applied,  in Section \ref{applications}, to  the fppf topos
and the Zariski topos of a scheme. Here we also discuss relations to 
the ``big'' fppf topos, which usually occurs in the literature.
In Section \ref{fppf-local} we introduce the notion of fppf-local rings
and establish their fundamental properties.
A construction of fppf-local rings depending on a given strictly local ring
and a limit ordinal is described in Section \ref{construction rings}.
This is used, in the final Section \ref{construction points}, to construct
topos-theoretical points for the fppf topos attached to a scheme.

\begin{remarkintro}
After this paper was submitted to the arXiv, I was kindly informed
by Shane Kelly that   related results   appear  in \cite{Gabber; Kelly 2015},
now a joint paper with Ofer Gabber.
\end{remarkintro}

\begin{acknowledgement}
I wish to thank the referee for noting some mistakes and  giving several suggestions, 
in particular for pointing out   that one has to distinguish between
the ``big'' fppf topos, which is usually used in the literature,
and the ``small'' fppf topos considered here.
\end{acknowledgement}

\section{Recollection: Universes, sites, topoi and locales}
\mylabel{universes sites topoi}

In this section, we  recall  some relevant foundational material about topos-theory from Grothendieck et.\ al.\ \cite{SGA 4a}.
Further very useful sources are Artin \cite{Artin 1962},  Johnstone \cite{Johnstone 1977},\cite{Johnstone 2002a},\cite{Johnstone 2002b}, 
Kashiwara and Shapira \cite{Kashiwara; Schapira 2006}, and the Stacks Project \cite{Stacks Project}.

Recall that a \emph{universe} is a nonempty set $\shU$ of sets satisfying four very natural
axioms, which we choose to state in the following form:  

\medskip
\begin{enumerate}
\item[(U  1)]
If $X\in\shU$ then $X\subset \shU$.
\smallskip
\item[(U 2)]
If $X \in\shU$ then $\left\{X \right\}\in\shU$.
\smallskip
\item [(U 3)]
If $X\in\shU$ then $\wp(X)\in\shU$.
\smallskip
\item [(U 4)]
If $I\in\shU$ and $X_i\in\shU$, $i\in I$ then $\bigcup_{i\in I}X_i\in\shU$.
\end{enumerate}

\medskip
In other words, $\shU$ is a nonempty transitive set of sets
that is stable under forming singletons, power sets, and unions indexed by 
$I\in\shU$. Roughly speaking, this ensures that   universes   are stable under the  set-theoretical operation
usually performed in practice.

If $X\in\shU$ is an arbitrary element, which  a fortiori exists because $\shU$ is nonempty, 
then the power set $\wp(X)$
and hence $\emptyset\in\wp(X)$ and the singleton $S=\left\{\emptyset\right\}$ are elements. 
In turn, the   set $I=\wp(S)$ of cardinality two is an element.
It follows by induction
that all finite ordinals    $0=\emptyset$ and   $n+1=n\cup\left\{n\right\}$ are elements as well.
Moreover, (U 2) and (U 4) ensure  that for each $X,Y\in \shU$, the set $\left\{X,Y\right\}$
is an element. The latter statement is the  form of the axiom (U 2) given in \cite{SGA 4a}, Expose I, Section 0.

Note that for a pair we have $(X,Y)\in\shU$ if and only if $X,Y\in\shU$, in light of Kuratowski's definition
of pairs $(X,Y)=\left\{\left\{X,Y\right\},\left\{Y\right\}\right\}$.
In particular, it follows that    groups or   topological spaces are elements of $\shU$ if and only if the underlying sets
are elements of $\shU$.
Note also that we adopt von Neumann's definition of \emph{ordinals} $\nu$ as sets of sets that are  transitive
(in the sense $\alpha\in\nu\Rightarrow \alpha\subset \nu$), so that the resulting order relation on $\nu$
 is a well-ordering (where  $\alpha\leq\beta$ means   $\alpha=\beta$ or $\alpha\in \beta$).
 In turn, each well-ordered set is order isomorphic to a unique ordinal.

Following Grothendieck, we assume
that any set $X$ is an element  of some universe $\shU$, which is an additional axiom of set theory.
Note that the intersection of universes is a universe, so there is always a unique smallest such universe.
The two axioms (U 3) and (U 4) enforce that the  cardinality $\aleph_\iota=\Card( \shU)$ of a universe
is \emph{strongly inaccessible}. In fact, the assumption that any set is 
contained in some universe is equivalent to the assumption that
any cardinal is majorized by some   strongly inaccessible cardinal.
A related notion was already mentioned by Felix Hausdorff  under the designation
``regul\"are Anfangszahlen mit Limesindex'', and I cannot resist from quoting  the original \cite{Hausdorff 1908}, page 131:
``Wenn es also regul\"are Anfangszahlen mit Limesindex gibt (und es ist uns bisher nicht
gelungen, in dieser Annahme einen Widerspruch zu entdecken), so ist die kleinste unter ihnen von einer so exorbitanten Gr\"o\ss e,
da\ss\ sie f\"ur die \"ublichen Zwecke der Mengenlehre kaum jemals in Betracht kommen wird.''

Given a universe $\shU$, a set $X$ is called a \emph{$\shU$-element} if $X\in\shU$.
We write $(\Set)_\shU$ for the category of all sets that are $\shU$-elements,
and likewise  denote by $(\Grp)_\shU, (\Sch)_\shU, (\Cat)_\shU $  
the categories  of all   groups, schemes, categories  that are $\shU$-elements.
The same notation is used for any other mathematical structure.
By  common abuse of notation, we sometimes drop the index, if there is no risk of confuction.
Note that a category $\calC$ is a $\shU$-element if and only if its
object set and all its hom sets have this property.
Given such a category, we denote by   $\PSh(\calC)$   the category of presheaves, that is, 
contravariant functors $\calC\ra(\Set)_\shU$. Given $X\in\calC$,
one writes the corresponding Yoneda functor as $h_X:\calC\ra(\Set)_\shU$, $Y\mapsto \Hom_\calC(Y,X)$.

A set $X$ is called \emph{$\shU$-small} if it is isomorphic to 
some element of $\shU$. The same locution is used for 
any mathematical structure, for example groups, rings, topological spaces,
schemes and categories.   
 If there is no risk of confusion, we simple use the term \emph{small} rather than $\shU$-small.
Note that a  category $\calC$ is a $\shU$-element or $\shU$-small if and only if
its object set and all its hom sets have the respective property.  
This has to be carfully distinguished from the following notion:
A category $\calC$ is called a \emph{$\shU$-category} if
the sets $\Hom_\calC(X,Y)$ are $\shU$-small for all objects $X,Y\in\calC$.
Clearly, this property is preserved by equivalences of categories.
It is also possible to define   a  Yoneda functor for
$\shU$-categories, and not only for categories that are $\shU$-elements,
by choosing bijections between $\Hom_\calC(Y,X)$ and elements of $\shU$.

Usually, the category $\PSh(\calC)$ is not $\shU$-small, even for $\calC\in\shU$:
Take, for example, the category of presheaves on a 
singleton space, which is equivalent to the category $(\Set)_\shU$.
Suppose its object set $\shU$ has the same cardinality $\kappa$
as an element $X\in\shU$. Since the power set $\wp(X)\in\shU$ has strictly larger
cardinality, and $\wp(X)\subset \shU$, we obtain $2^\kappa\leq \kappa$, in 
 contradiction  to cardinal arithmetic. 

Let $\calC$ be a category.
A \emph{Grothendieck topology} on $\calC$ is 
a collection  $J(X)$ of sieves for each object $X\in\calC$, satisfying certain
axioms. We do not bother to reproduce the axioms, and refer for details to   \cite{SGA 4a}, Expose II, Section 1.
Recall that a \emph{sieve} on $X$  is a full subcategory $\calS\subset\calC/X$
with the property that  $A\in\calS$ for each morphism $A\ra B$ in $\calC/X$ with $B\in\calS$.
Usually,   the covering sieves of a Grothendieck topology are specified with the help of \emph{pretoplogies},
which is a collection $\Cov(X)$ of tuples $(X_\alpha\ra X)_{\alpha\in I}$
of morphisms $X_\alpha\ra X$ for each object $X\in\calC$ satisfying similar axioms.
These tuples are  referred to  as \emph{coverings} of $X$, and the
induced   Grothendieck topology is the finest one for which the
coverings families generate coverings sieves.

A category $\calC$ endowed with a Grothendieck topology is called a \emph{site}.
To proceed, choose a universe with $\calC\in\shU$.
Then  we have the full subcategory $\Sh(\calC)\subset\PSh(\calC)$
of sheaves, that is, contravariant functors $\calC\ra(\Set)_\shU$ satisfying
the sheaf axioms with respect to the covering sieves or covering families.
A \emph{$\shU$-site} is a $\shU$-category $\calC$ endowed with a Grothendieck topology,
so that there is a full $\shU$-small subcategory $\calD\subset\calC$ so that
for each object $X\in\calC$, there is a covering family $(X_\alpha\ra X)_\alpha$ with $X_\alpha\in\calD$.
This condition ensures that the category $\Sh(\calC)$ of $\shU$-sheaves remains a $\shU$-category.

A    \emph{$\shU$-topos} is a $\shU$-category $\shE$ that is equivalent
to the category $\Sh(\calC)$ for some site $\calC\in\shU$.
The central result on topoi is \emph{Giraud's Characterization} (\cite{SGA 4a}, Expose IV, Theorem 1.2).
Roughly speaking, it makes the following three assertions:
First, it says that one may choose the site $\calC$  so that it
contains all inverse limits, and that its Grothendieck topology is \emph{subcanonical},
which means that all Yoneda functors $h_X$, $X\in\calC$ satisfy the sheaf axioms.
Second, it singles out the  topoi among the $\shU$-categories
in terms of   purely categorical properties of $\shE$, referring to objects and arrows
rather than to coverings.
Third, it characterizes the topoi among the $\shU$-sites
using the \emph{canonical topology} on $\shE$, which is 
the finest topology on $\shE$ that turns all Yoneda functors $h_F$, $F\in\shE$ into sheaves.

A \emph{continuous map} $\epsilon:\shE\ra\shE'$ between $\shU$-topoi 
is a triple $\epsilon=(\epsilon_*,\epsilon^{-1},\varphi)$,
where $\epsilon_*:\shE\ra\shE'$  and $\epsilon^{-1}:\shE'\ra\shE$ are
adjoint functors, and $\varphi$ is the adjunction isomorphism.
Here $\epsilon^{-1}$ is left adjoint and called the \emph{preimage functor}, and $\epsilon_*$ is right adjoint
and called the \emph{direct image functor}.
Moreover, one demands that the  preimage functor  $\epsilon^{-1}$ is left exact, that
is, commutes with finite inverse limits.
Up to isomorphism, the continuous map $\epsilon$ is determined by
either of the preimage functor $\epsilon^{-1}$ and the direct image functor $\epsilon_*$.
The set of all continuous maps $\Hom(\shE,\shE')$ is itself a category,
the morphism being the compatible natural transformations between the direct image
and preimage functors.

The $\shU$-valued sheaves on a topological space $X\in\shU$ form a a topos $\shE=\Sh(X)$.
In particular, the category $(\Set)_\shU$ can be identified with 
the category $\Sh(S)$ for the singleton space $S=\left\{\star\right\}$.
In light of this, a \emph{point in the sense of topos-theory} of a $\shU$-topos $\shE$ is a continuous map of topoi $P:(\Set)_\shU\ra \shE$.
We denote by $\Points(\shE)$ the category of points,
and by $|\shE|$ the set of isomorphism classes $[P]$ of points.
This set is endowed with a natural topology:
Choose a terminal object $e\in\shE$. Given  a subobject $U\subset e$,
we formally write  $U\cap |\shE|\subset |\shE|$ for the set of
isomorphism classes of points $P$ with $P^{-1}(U)\neq\emptyset$, and declare
it as open. This   indeed constitutes a topology on the set $|\shE|$.

Recall that a nonempty ordered set $L$ is called a  \emph{locale}
if the following axioms hold: 

\medskip
\begin{enumerate}
\item[(LC 1)]
For all pairs $U,V\in L$, the infimum $U\wedge V\in L$ exists, that
is, the largest element that is smaller or equal than both $U,V$.
\smallskip
\item[(LC 2)]
For each family $U_i\in L$, $i\in I$,
the supremum  $\bigvee_{\alpha\in I} U_i\in L$ exists, that
is, the least element that is larger or equal than all $U_i$.
\smallskip
\item [(LC 3)]
The   distributive law holds, which means $U\wedge(\bigvee_{i\in I} V_i)=\bigvee_{i\in I}(U\wedge V_i)$.
\end{enumerate}

\medskip
The ordered set  $L=\shT$ of open subsets $U\subset X$ of a topological space $(X,\shT)$ is the paramount
example for locales.  One should regard 
locales as abstractions of topological spaces, where one  drops the underlying set and merely keeps the topology.
In light of this, one defines a \emph{continuous map} $f:L\ra L'$ between locales as
a monotonous map $f^{-1}:L'\ra L$   that 
respects finite infima and arbitrary suprema. Note the reversal of arrows. Here the notation $f^{-1}$
is purely formal, and does not indicate that $f$ is bijective. 

Each locale comes with a Grothendieck topology, and thus can be regarded as site:
The covering families $(U_\alpha\ra V)_{\alpha\in I}$ are
those with $V=\bigvee_\alpha U_\alpha$. In turn, we have the $\shU$-topos $\Sh(L)$
of sheaves on the locale $L\in\shU$.
Conversely, for each  $\shU$-topos $\shE$ we have 
a locale $\Loc(\shE)$, which is the ordered set of 
subojects $U\subset e$ of a fixed terminal object $e\in\shE$.
Up to canonical isomorphism, it does not depend on the choice of the terminal object.

\section{Topoi with   same locales}
\mylabel{topoi same locales}

In this section, we establish some   facts on continuous maps between certain
topoi, which   occur in algebraic geometry when various Grothendieck topologies are involved.
In order to achieve flexibility and facilitate application, we work in the following axiomatic set-up: 

Throughout,   fix a universe $\shU$. 
Let $\calC_f$ and $\calC_z$ be two    categories that are $\shU$-elements, in which  
terminal objects exist. Furthermore, suppose these categories are equipped
with a pretopology of coverings, such that we regard them as sites.
We now suppose that we have   adjoint functors
$$
u:\calC_f\lra\calC_z\quadand v:\calC_z\lra\calC_f,
$$
where $u$ is the left adjoint and  $v$ is the right adjoint.  
Let us write the objects of $\calC_f$ \emph{formally} as pairs $(U,p)$, and the objects of
$\calC_z$ by ordinary letters $V$, and the functors as
$$
u(U,p)=p(U)\quadand v(V)=(V,i).
$$
Note that the adjunction, which by abuse of notation is regarded as an identification,   takes the form
\begin{equation}
\label{adjoint functors}
\Hom_{\calC_z}(p(U),V) = \Hom_{\calC_f}((U,p),(V,i)).
\end{equation}

Let me emphasize  that all \emph{this notation is purely formal}, 
but   based on   geometric intuition.
The guiding example, which one should have in mind, is that
$\calC_z$ comprises   open subsets of a scheme $X$,
and that $\calC_f$ consists of certain   flat $X$-schemes $(U,p)$, where $p:U\ra X$ is the structure
morphism that is assumed to be universally open. The functor $u$ takes such an $X$-scheme to its image $p(U)\subset X$,
whereas the functor $v$ turns the open subset $V\subset X$ into an $X$-scheme $(V,i)$,
where $i:V\ra X$ is the inclusion morphism. Note that, by abuse of notation, we usually write
$i:V\ra X$ and not the more precise $i_V:V\ra X$.
Of course, the indices in $\calC_f$ and $\calC_z$
refer to ``flat'' and ``Zariski'', respectively.

We now demand the following four conditions (TL 1) -- (TL 4), which conform with    geometric intuition:

\medskip
\begin{enumerate}
\item[(TL 1)]
The composite functor $u\circ v$ is isomorphic to the identity on $\calC_z$. 
\smallskip
\item[(TL 2)]
For each covering family $(V_\lambda\ra V)_\lambda$ in the site $\calC_z$,
the induced family $((V_\lambda,i_\lambda)\ra (V,i))_\lambda$   in the site $\calC_f$ is covering.
\smallskip
\item [(TL 3)]
For each family $(U_\lambda,p_\lambda)_\lambda$  of objects in $\calC_f$, there is a subobject $V$
of the terminal object in $\calC_z$ and factorizations $p_\lambda(U_\lambda)\ra V$ 
so that the   family $(p_\lambda(U_\lambda)\ra V)_\lambda$ is a covering.
\smallskip
\item[(TL 4)]
For each object $(U,p)$ in $\calC_f$, the representable presheaf $h_{(U,p)}$ is a sheaf.
\end{enumerate}

\medskip
Let me make the following remarks:
The functor $v:\calC_z\ra\calC_f$, being a right adjoint,  commutes with inverse limits.
Moreover, condition (TL 1) ensures that $v$ is   faithful, which allows us to  make the
identification $i(V)=V$.
By our overall assumption on  the categories $\calC_f$ and $\calC_z$, the   terminal object  appearing in  condition (TL 3) does exists.
Finally, condition (TL 4) can be rephrased as that the Grothendieck  topology on $\calC_f$ is
finer than the \emph{canonical topology}, which is the finest topology for which every representable presheaf
satisfies the sheaf axioms. One also says that the Grothendieck topology on $\calC_f$ is \emph{subcanonical}.

Furthermore, I want to point out that we do not assume that our categories $\calC_f,\calC_z$ all finite
inverse limits are representable. However, it is part of the definition for
\emph{pretopologies} that for any covering family $(U_\lambda\ra U)_\lambda$,
the members $U_\lambda\ra U$ are \emph{base-changeable} (``quarrable''  in \cite{SGA 4a},   Expose II, Definition 1.3),
that is, for every other morphism $U'\ra U$, the fiber product $U_\lambda\times_UU'$ does exist.

\begin{proposition}
\mylabel{cocontinuous and continuous}
The functor $u:\calC_f\ra\calC_z$ is cocontinuous, and the adjoint functor $v:\calC_z\ra\calC_f$
is continuous.
\end{proposition}

\proof
For the precise definition of \emph{continuous} and \emph{cocontinuous functors} between sites, we refer
to \cite{SGA 4a}, Expose III. 
The two assertions are equivalent, according to  loc.\ cit.\ Proposition 2.5, because the functors $u$ and $v$ are adjoint.
To check that $u$ is cocontinuous, let $(U,p)\in\calC_f$, and $(V_\lambda\ra p(U))_\lambda$ be a covering family in $\calC_z$.
By condition (TL 2), the  induced family  $((V_\lambda,i_\lambda)\ra (p(U),i))_\lambda$ is a covering family in $\calC_f$.
Form the pull-back
$$
\begin{CD}
(U_\lambda,p_\lambda) @>>> (V_\lambda,i_\lambda)\\
@VVV @VVV\\
(U,p) @>>> (p(U),i)
\end{CD}
$$
in $\calC_f$, which exists because members of covering families are base-changeable.
By the axioms for covering families,  $((U_\lambda,p_\lambda)\ra (U,p))_\lambda$ remains a covering family.
The preceding diagram, together with the adjunction,
shows that  the induced maps $p_\lambda(U_\lambda)\ra p(U)$ factor over $V_\lambda\ra p(U)$.
 If follows that $u$ is cocontinuous (\cite{SGA 4a}, Expose III, Definition 2.1).
\qed

\medskip
Now let $\shE_f=\Sh(\calC_f)$ and $\shE_z=\Sh(\calC_z)$ be the $\shU$-topoi of sheaves on $\calC_f$ and $\calC_z$, respectively.
We refer to the sheaves on $\calC_f$ as \emph{f-sheaves},
and to the sheaves on  $\calC_z$ as \emph{z-sheaves}.
According to \cite{SGA  4a}, Expose IV, Section 4.7, the cocontinuous functor $u:\calC_f\ra\calC_z$ induces a continuous map of topoi
$$
\epsilon=(\epsilon_*,\epsilon^{-1},\varphi):\shE_f\lra\shE_z.
$$
Let me make this explicit: The direct image functor $\epsilon_*$ sends an f-sheaf $F$ to the z-sheaf $\epsilon_*(F)$ defined by
\begin{equation}
\label{direct image}
\Gamma(V,\epsilon_*(F))=\Gamma((V,i),F).
\end{equation}
Note that, in general, this would be
an inverse limit of local sections, indexed by the category of pairs $((U,p),\psi)$,
where $\psi:p(U)\ra V$ is a morphism in $\calC_z$. In our situation,
such an inverse limit is not necessary, because the index category has a terminal object
$((V,i),\psi)$, where $\psi:i(V)\ra V$ is the canonical isomorphism coming from condition (TL 1).

The inverse image functor $\epsilon^{-1}$ sends a z-sheaf $G$ to the f-sheaf $\epsilon^{-1}(G)$,
defined by
\begin{equation}
\label{inverse image}
\Gamma((U,p),\epsilon^{-1}(G))=\Gamma(p(U),G).
\end{equation}
Note that, in general, this would give merely a presheaf, and   sheafification is necessary.
In our situation, however, the presheaf is already a sheaf, thanks to condition (TL 2). 

The adjunction map $\varphi$ between $\epsilon^{-1}$ and $\epsilon_*$ is determined by   natural transformations
$$
G\lra \epsilon_*\epsilon^{-1}(G)\quadand \epsilon^{-1}\epsilon_*(F)\lra F.
$$
Here the former comes from identity maps
$$
\Gamma(V,G)\stackrel{\id}{\lra}\Gamma(i(V),G)=\Gamma((V,i),\epsilon^{-1}(G))=\Gamma(V,\epsilon_*\epsilon^{-1}(G)).
$$
The latter is the given by   restriction maps
$$
\Gamma((U,p),\epsilon^{-1}\epsilon_*(F))=\Gamma(p(U),\epsilon_*(F))=\Gamma((p(U),i),F)\stackrel{\res}{\lra}\Gamma((U,p),F).
$$
For later use, we now establish a technical fact:

\begin{lemma}
\mylabel{technical fact}
For each object $V$ in $\calC_z$, the presheaf $h_V$ is a sheaf. Furthermore, we have $\epsilon^{-1}(h_V)=h_{(V,i)}$.
\end{lemma}

\proof
According to condition (TL 4), the presheaf $h_{(V,i)}$ on $\calC_z$ is a sheaf. If the follows from
condition (TL 2) and the fact that the functor $v:\calC_z\ra\calC_f$ is   faithful that the presheaf
$h_V$ on $\calC_z$ satisfies the sheaf axioms. Finally, for each object $(U,p)$ of $\calC_f$, we have 
$$
\Gamma((U,p),\epsilon^{-1}(h_V))=\Gamma(p(U),h_V)  = \Hom_{\calC_z}(p(U),V),
$$
where the first equation comes form (\ref{inverse image}), and the second equation stems
from the Yoneda Lemma. Similarly, we have
$$
\Gamma((U,p),h_{(V,i)}) = \Hom_{\calC_f}((U,p),(V,i)).
$$
Using the adjointness of the functors (\ref{adjoint functors}), we infer that $\epsilon^{-1}(h_V)=h_{(V,i)}$.
\qed

\medskip
The continuous map of topoi $\epsilon:\shE_f\ra\shE_z$ induces a functor  $\epsilon:\Points(\shE_f)\ra\Points(\shE_z)$ on the category
of points. In turn, we get a  continuous map of topological spaces
$\epsilon: |\shE_f|\lra |\shE_z|$.
To understand this map, we first look at the induced continuous map of locales
\begin{equation}
\label{morphism locale}
\epsilon:\Loc(\shE_f)\lra\Loc(\shE_z).
\end{equation}
Recall that these locales comprise the 
ordered sets of subobjects of the chosen terminal objects.
Let us write  $X\in \calC_z$ for the terminal object.
It follows that $e_z=h_{X}$ is a terminal object in the topos $\shE_z$, and $\Loc(\shE_z)$ is the ordered
set of subobjects $G\subset e_z$.
Being right adjoint, the functor $v:\calC_z\ra\calC_f$ respects inverse limits, whence $(X,i)\in\calC_f$ is a terminal object.
In turn, $e_f=h_{(X,i)}$ is the terminal object in the topos $\shE_f$, and $\Loc(\shE_f)$ is the ordered set of subobjects $F\subset e_f$.
The continuous map of locales (\ref{morphism locale})  is  just the monotonous map
\begin{equation}
\label{map  locale}
\Loc(\shE_z)\lra\Loc(\shE_f),\quad G\longmapsto \epsilon^{-1}(G),
\end{equation}
which, by definition of the hom-sets in the category of locales, goes in the reverse direction.

\begin{theorem}
\mylabel{locale isomorphic}
The continuous map of locales $\epsilon:\Loc(\shE_f)\ra\Loc(\shE_z)$ is a homeomorphism,
that is, the monotonous map (\ref{map  locale}) is bijective.
\end{theorem}

\proof
The argument is analogous to \cite{SGA 4b}, Expose VIII, Proposition 6.1.
To see that the monotonous map is injective, suppose we have two subobjects $G,G'\subset e_z$
with $\epsilon^{-1}(G)=\epsilon^{-1}(G')$ as subobjects in $e_f=\epsilon^{-1}(e_z)$.
For each $V\in\calC_z$, we then have
$$
\Gamma(V,G)=\Gamma((V,i),\epsilon^{-1}(G))=\Gamma((V,i),\epsilon^{-1}(G'))=\Gamma(V,G'),
$$
where the outer identifications come from (\ref{inverse image}). Whence $G=G'$.

For surjectivity, let $F\subset e_f$ be a subobject  in $\calC_f$.
First note that, for each object $(U,p)$ in $\calC_f$ over which $F$ has a local section,
the set $\Gamma((U,p),F)$ must be a singleton, because $F$ is a subobject of the terminal object.
Moreover, if $(U',p')\ra (U,p)$ is a morphism, then $\Gamma((U',p'),F)$ stays a singleton.

Now consider the family $(U_\lambda,p_\lambda)$ of all objects in $\calC_f$ over which $F$ has
a local section. Using condition (TL 3), there is a subobject $V\subset X$ of the terminal object and morphisms
$p_\lambda(U_\lambda)\ra V$ so that the induced family $((U_\lambda,p_\lambda)\ra (V,i))_\lambda$ is covering.
Consider the fiber products
$$
\begin{CD}
(U_{\lambda\mu},p_{\lambda\mu}) @>>> (U_\mu,p_\mu)\\
@VVV @VVV\\
(U_\lambda,p_\lambda) @>>> (V,i),
\end{CD}
$$
which exists because members of covering families are base-changeable.
The sheaf axioms give a short exact sequence
$$
\Gamma((V,i),F)\lra\prod_{\lambda}\Gamma((U_\lambda,p_\lambda),F)\lra\prod_{\lambda,\mu}\Gamma((U_{\lambda\mu},p_{\lambda\mu}),F),
$$
where the terms in the middle and the right are singletons. If follows that the term on the left is a singleton.
The Yoneda Lemma yields a morphism of presheaves $h_{(V,i)}\ra F$. Note that the presheaf $h_{(V,i)}$ is actually a sheaf,
according to condition (TL 4).  Moreover, this morphism is actually an isomorphism:
Let $(U,p)$ be an arbitrary object of $\calC_f$. Suppose there is a morphism $(U,p)\ra (V,i)$.
Since $V$ is a  subobject of the terminal object, so is $(V,i)$, by the adjunction (\ref{adjoint functors}). It follows that the term on the left in
\begin{equation}
\label{map of local sections}
\Gamma((U,p),h_{(V,i)})\lra\Gamma((U,p),F)
\end{equation}
is a singleton, whence the map is bijective. Finally, suppose there is no morphism $(U,p)\ra (V,i)$,
such that there is also no morphism $p(U)\ra V$. By the very definition of $V$, this means $\Gamma((U,p),F)=\emptyset$.
Again, the map (\ref{map of local sections}) is bijective. We conclude that $h_{(V,i)}\ra F$ is an isomorphism.
Using Lemma \ref{technical fact}, we infer that $F=\epsilon^{-1}(h_V)$. Since $V$ is a subobject of
the terminal object in $\calC_z$, the sheaf $h_V$ must be a subobject of the terminal object $e_z\in\shE_z$,
which concludes the proof.
\qed

\medskip
We finally come to the induced continuous map $\epsilon:|\shE_f|\ra|\shE_z|$
of topological spaces.  Recall that the \emph{chaotic topology} on a set 
has as sole open subsets the whole set and the empty set.
A topological space $X$ is called \emph{sober} if each irreducible closed subset has a unique
generic point. Each space $X$ comes with the \emph{sobrification} $X\ra X_\sob$, which
is universal   with respect to continuous maps into sober spaces,
compare \cite{EGA I}, Chapter 0, Section 2.9.

\begin{corollary}
\mylabel{fiber chaotic}
Each fiber of the map $\epsilon:|\shE_f|\ra|\shE_z|$
carries the chaotic topology, and 
the induced map     of sober spaces
is an embedding $|\shE_f|_\sob\subset |\shE_z|_\sob$. The latter is an equality
provided that the map $\epsilon:|\shE_f|\ra|\shE_z|$ admits a section.
\end{corollary}

\proof
According to the theorem, each open subset in $|\shE_f|$ is the preimage of
an open subset in $|\shE_z|$. The statement is thus a special case of the following lemma.
\qed

\begin{lemma}
\mylabel{bijective locales}
Let $f:X\ra Y$ be a continuous map between topological spaces. Suppose that
each open subset  in $X$ is the preimages of an open subset  in $Y$.
Then all fibers of $f$ carry the chaotic topology, and the induced
map of sober space is an embedding $X_\sob\subset Y_\sob$.
The latter is an equality provided that $f$ admits a section.
\end{lemma}

\proof
Given $y\in Y$, and let $U=f^{-1}(V)$ be an open subset in $X$.
Then $U\cap f^{-1}(y)$ is either empty or the hole fiber. In turn, the fiber $f^{-1}(y)$ carries
the chaotic topology. 
Likewise, one sees that $f:X\ra Y$ is a closed map.
Recall that one may view the sobrification $Y_\sob$ as the space of closed irreducible subsets in $Y$.
Let $Z\subset Y$ be such a subset. Then $f(Z)\subset X$ is a closed irreducible subset as well.
Given $x\in f^{-1}f(Z)$ and any open neighborhood $x\in U=f^{-1}(V)$, we see that 
$f(f^{-1}(V)\cap Z)=V\cap f(Z)$ is nonempty. Whence $x$ is in the closure of $Z$, which means $x\in Z$.
We conclude that $Z\subset X$ equals the preimage of the closed irreducible subset $f(Z)$.
In turn, we may regard the induced map on sobrification   
as an inclusion $X_\sob\subset Y_\sob$. The space $X_\sob$ carries the subspace topology,
again by our assumption 
on the open sets.

One easily sees  that any set-theoretical section $s:Y\ra X$ for $f$ must be continuous.
It thus induces a right inverse for the inclusion $X_\sob\subset Y_\sob$, which then must
must be an equality.                                                                         
\qed

\section{Applications to algebraic geometry}
\mylabel{applications}

We now   apply the abstract results of the preceding section to some concrete Grothendieck topologies
in algebraic geometry.
Let $X$ be a scheme, and fix a universe with $X\in\shU$. In what
follows, all schemes are tacitly assumed to be $\shU$-elements.
We denote by $(\Zar/X)$ the locale given by the  ordered set of  open subschemes $V\subset X$,
regarded as a site in the usual way, and write  $X_\Zar$ for the ensuing $\shU$-topos of sheaves on $X$.

Let us denote by
$(\fppf/X)$
the category of $X$-schemes $(U,p)$, where the structure morphism $p:U\ra X$  is
 flat and locally of finite presentation, following the convention of \cite{Stacks Project}.
The hom sets in this category are formed by arbitrary $X$-morphisms.
Note that any such morphism is automatically  locally of finite presentation, by \cite{EGA I}, Proposition 6.2.6,
but not necessarily flat.
If $(U,p),(V,q)$ are two objects and $U\ra V$ is a flat $X$-morphism,
than for any other  object $(V',q')$ and any $X$-morphism $V'\ra V$,
the usual fiber product of schemes $U\times_VV'$ yields an object and whence a fiber product in $(\fppf/X)$.
 It is not clear to me to what extend
other fiber products in $(\fppf/X) $ exists, which may differ from the ususal fiber products in $(\Sch/X)$.

Our category is equipped with the pretopology of \emph{fppf coverings}
$((U_\alpha,p_\alpha)\ra (U,p))_{\alpha\in I}$, where each $U_\alpha\ra U$ is
flat, and the induced map $\coprod_\alpha U_\alpha\ra U$ is surjective.
We regard $(\fppf/X)$ as a site, and denote by  $X_\fppf$ the resulting $\shU$-topos
of sheaves.  Note that the category $(\fppf/X)$ usually contains hom sets of cardinality $\geq 2$,
in contrast to the Zariski   site.

Clearly,  the categories $(\Zar/X)$ and $(\fppf/X)$ have terminal objects.
Consider the functors
$$
u:(\fppf/X)\lra(\Zar/X),\quad (U,p)\longmapsto p(U)
$$
and
$$
v:(\Zar/X)\lra(\fppf/X),\quad V\longmapsto (V,i)
$$
where $i:V\ra X$ denotes the inclusion morphism of an open subscheme.
These are well-defined, because any flat morphism locally of finite presentation is universally open
(\cite{EGA IVb}, Theorem 2.4.6), 
and any open embedding is a fortiori flat and locally of finite
presentation. Note that, for schemes that fail to be quasiseparated, there are open 
subschemes whose inclusion morphism is not quasicompact, and in particular not of finite presentation.
Nevertheless, they are locally of finite presentation.

\begin{proposition}
\mylabel{fppf axioms}
The functor $u$ is left adjoint to $v$, and this pair of adjoint functors
satisfy the conditions (TL 1) -- (TL 4) of Section \ref{topoi same locales}.
\end{proposition}

\proof
The adjointness follows from the universal property of schematic images of open morphisms.
The first two conditions (TL 1) and (TL 2) are trivial. To see (TL 3),
let $p_\lambda:U_\lambda\ra X$ be flat and locally of finite presentation.
The the image is open, and the open subscheme $V=\bigcup_\lambda p_\lambda(U_\lambda)$ of the terminal object $X\in(\Zar/X)$ is
covered by the $p_\lambda(U_\lambda)$.
Finally, (TL 4) holds by \cite{SGA 1}, Expose VIII, Theorem 5.2. 
\qed

\medskip
In turn, the functor $u:(\fppf/X)\ra(\Zar/X)$ is cocontinuous and induces a   morphism of topoi $\epsilon:X_\fppf\ra X_\Zar$.
Applying Theorem \ref{locale isomorphic} and its Corollary we get:

\begin{proposition}
\mylabel{fppf homeomorphisms}
The induced continuous map   $\epsilon:\Loc(X_\fppf)\ra\Loc(X_\Zar)$ of locales is a homeomorphism,
and the induced map of sober spaces is an embedding $ |X_\fppf|_\sob\subset |X_\Zar|_\sob=X$.
\end{proposition}

In Section \ref{construction points}, we shall construct a section for the map $|X_\fppf|\ra|X_\Zar|$,
whence we actually have an equality $|X_\fppf|_\sob=|X_\Zar|_\sob=X$.

This approach carries over to  analogous topoi defined with the  \'etale topology, the Nisnevich 
topology \cite{Nisnevich 1989},
and the syntomic topology \cite{Fontaine; Messing 1987}. With the obvious notation, we thus get  canonical identifications
of locales and sober topological  spaces. We leave the details to the reader. 
Note, however, that Theorem \ref{locale isomorphic} does not apply to the fpqc topology.
This is because   flat and quasicompact morphisms are not necessarily open (for example the one induced by 
the faithfully flat ring extension $\ZZ\subset\ZZ\times\QQ$,
compare \cite{EGA IVb}, Remark 2.4.8). Thus  we apparently have no functor $v$ from the fpqc site to the Zariski site.
However, one may apply it to the site $(\fpuo/X)$ of \emph{ flat and universally open morphisms},
which was considered by Romagny \cite{Romagny 2003}.
It would be interesting to understand the relation between the fpqc topos and the fpuo topos.

One disadvantage for the site $(\fppf/X)$ and the enusing topos $X_\fppf$ is that
it is apparently not functorial with respect to $X$, by similiar reasons
as for the lisse-\'etale topos. Therfore, in the literature one usually considers the
category $(\Sch/X)$ of all $X$-schemes contained in a fixed universe, 
and endows it with the fppf topology. The covering families
are the $(U_\alpha\ra Y)_\alpha$, where   $U_\alpha\ra Y$ are flat  $X$-morphisms that
are locally of finite presentation, and $\amalg U_\alpha\ra Y$ is surjective.
Let me refer to sheaves on this site $(\Sch/X)$ as 
 \emph{big fppf sheaves}, whereas we now call sheaves
on the site $(\fppf/X)$ \emph{small fppf sheaves}.
Likewise, we call the resulting topos 
$X_\fppf^\bg$ the \emph{big fppf topos}, whereas we refer to $X_\fppf$ as the  \emph{small fppf topos}.

Given a big fppf sheaf $F$ over $X$, we obtain by forgetting superfluous local sections and restriction maps
a small fppf sheaf $F|Y$ for each $X$-scheme $Y$. The resulting functor
\begin{equation}
\label{forgetful}
X_\fppf^\bg\lra Y_\fppf,\quad F\longmapsto F|Y
\end{equation}
commutes with direct and inverse limits, because the sheafification functor $F\mapsto F^{++}$ involves
over a fixed $Y$ only fppf coverings and their fiber products, which are the same in $(\Sch/X)$ and $(\fppf/Y)$.
In turn, we obtain a functor
$$
\Points(Y_\fppf)\lra\Points(X_\fppf^\bg),\quad P\longmapsto P^\bg,
$$
where $P^\bg$ is given  by the fiber functor $F_{P^\bg}=(F|Y)_P$.
Likewise, we get a continuous map of locales: Fix a terminal object $e\in X_\fppf^\bg$.
Then $e|Y$ are terminal objects, and for each subobject $G\subset e$
we get subojects $G|Y\subset e|Y$. Since the forgetful functor (\ref{forgetful}) commutes
with direct and inverse limits, the monotonous map $G\mapsto G|Y$ constitutes a continuous map of locales
\begin{equation}
\label{locales}
\Loc(Y_\fppf)\lra \Loc(X_\fppf^\bg). 
\end{equation}
In turn, the canonical map of topological spaces $|Y_\fppf|\ra|X_\fppf^\bg|$ is continuous as well.

\section{fppf-local rings}
\mylabel{fppf-local}

We now define a class of local rings that generalizes 
the notion of algebraically closed fields. Such rings will appear   as
stalks of the structure sheaf $\O_{X_\fppf}$ at topos-theoretical points.

\begin{definition}
\mylabel{fppf-local rings}
A ring $R$ is called  \emph{fppf-local} if it is local,
and every fppf homomorphism $R\ra B$ admits a retraction.
In other words, the corresponding morphism of schemes
$\Spec(B)\ra\Spec(R)$ admits a section.
\end{definition}

Recall that a ring $R$ is called \emph{totally integrally closed} 
if for any ring homomorphism $B\ra R$  
and any integral extension $B\subset B'$, there is homomorphism $B'\ra R$ making the
diagram 
$$
\xymatrix{
			& \Spec(B')\ar[d]\\
\Spec(R)\ar[r]\ar[ur]	& \Spec(B)
}
$$
commutative. This was introduced by Enochs \cite{Enochs 1968}, and further
analyzed by Hochster \cite{Hochster 1970}. One also says that $R$ is a \emph{TIC ring}.
Note that such rings are necessarily reduced (see \cite{Enochs 1968}, Theorem 1 and also \cite{Borho; Weber 1971}).

Let us call a ring $R$  \emph{absolutely integrally closed}
if each monic polynomial $f\in R[T]$ has a root in $R$.
These rings are also called \emph{AIC rings}.
Note that an integral domain $R$ is AIC  if and only its it is normal and its field of fraction
is algebraically closed, and this holds if and only if $R$ is TIC (\cite{Enochs 1968}, Proposition 3).
Throughout, we call a ring  $R$ \emph{integral} if it is an integral domain, that
is, a subring of a field.

\begin{proposition}
\mylabel{TIC fppf-local}
Any fppf-local ring is AIC. Moreover, for  local rings $R$, the following three conditions are equivalent:
\begin{enumerate}
\item $R$ is TIC.
\item $R$ is AIC and integral.
\item $R$ is fppf-local and integral.
\end{enumerate}
\end{proposition}

\proof
Suppose that $R$ is fppf-local, and let $f\in R[T]$ be a monic polynomial.
The fppf algebra $R[T]/(f)$ contains a root of $f$, and also admits an $R$-algebra homomorphism
to $R$. Hence $R$ itself contains a root. Consequently, $R$ is AIC.
This also shows that implication (iii)$\Rightarrow$(ii).
According to \cite{Hochster 1970}, Proposition 7, a local ring   that is  TIC must be integral.
The equivalence of (i) and (ii) follows from \cite{Enochs 1968}, Proposition 3.
Now suppose that $R$ is TIC. Using the TIC condition with $B=R$ and a finite fppf algebra $B'$,
we infer that $R$ is fppf-local. This gives the implication (i)$\Rightarrow$(iii).
\qed

\medskip
We shall see in Section \ref{construction rings} that there are fppf-local rings that are not integral.
I do not know whether there are local  AIC rings that are not fppf-local.
Neither do I know whether nonzero homomorphic images of fppf-local rings remain fppf-local.
We have the following partial results in this direction:

\begin{proposition}
\mylabel{fppf residue}
Let $R$ be an fppf-local ring. For every prime ideal $\primid\subset R$, the   domain $R/\primid$
is fppf-local, and the residue field $\kappa(\primid)$ is algebraically closed.
For every ideal $\ideala\subset R$,
the  residue class ring $R/\ideala$ is  AIC.
\end{proposition}

\proof
Let $\bar{f}\in R/\ideala[T]$ be a monic polynomial. Lift it to a monic polynomial $f\in R[T]$.
The fppf $R$-algebra $B=R[T]/(f)$ contains a root of $f$, and
admits an $R$-algebra homomorphism to $R$. Whence there is root $a\in R$,
whose residue class is a root of $\bar{f}$ in $R/\ideala$, such that the latter is AIC.
If $\ideala=\primid$ is prime, then the AIC domain $R/\primid$ is fppf-local
by Proposition \ref{TIC fppf-local}.
Any localization of the domain $R/\primid$ stays AIC. In particular, its field of fractions $\kappa(\primid)$
is algebraically closed.
\qed

\begin{proposition}
\mylabel{cotangent}
Let $R$ be an fppf-local ring, and $\primid\subset R$ be a prime ideal.
Then $\primid^n=\primid$ for every integer $n\geq 1$.
\end{proposition}

\proof
Let $a\in\primid$ be some element, and consider the monic polynomial $f=T^n-a\in R[T]$.
The fppf $R$-algebra $R[T]/(f)$ contains a root for $f$, and admits an $R$-algebra homomorphism
to $R$, whence there is an element $b\in R$ with  $b^n=a$. Since $\primid$ is
prime, we must have $b\in\primid$.
\qed

\medskip
In particular, all cotangent spaces $\primid/\primid^2$ of an fppf-local ring $R$ vanish,
and its formal completion $\hat{R}=\invlim_n R/\maxid^n$ coincides with the residue field $\kappa=R/\maxid$.
This yields the following:

\begin{corollary}
\mylabel{fppf-local noetherian}
A ring is fppf-local and noetherian if and only if it is an algebraically closed field.
\end{corollary}

\proof
The condition is   sufficient by Proposition \ref{TIC fppf-local}, or more directly by  Hilbert's Nullstellensatz. 
Conversely, suppose that $R$ is an fppf-local noetherian ring, with residue field $k=R/\maxid$.
By the Proposition,  $\maxid\otimes k=\maxid/\maxid^2=0$. The Nakayama Lemma ensures $\maxid=0$,
whence $R=k$ is a field. This field is algebraically closed by Proposition \ref{TIC fppf-local}.
\qed

\medskip
It follows that finite flat algebras over fppf-local rings may fail to be fppf-local:
Take $R=k$ an algebraically closed field and $A=k[\epsilon]$ the ring of dual numbers,
where $\epsilon^2=0$.  

We now can state an amazing property of fppf-local rings:

\begin{theorem}
\mylabel{amazing property}
Let $R$ be an fppf-local ring. For every prime ideal $\primid\subset R$,
the local ring $R_\primid$ is strictly local, with algebraically closed residue field.
\end{theorem}

\proof
The residue field $\kappa(\primid)$ is algebraically closed by Proposition \ref{fppf  residue}.
It remains to check that the local rings $R_\primid$ are henselian.
Let  
$\idealq\subset R$ be a minimal prime ideal contained in $\primid$. 
Then the local domain $(R/\idealq)_\primid$ is AIC according to  Proposition \ref{fppf residue},
whence henselian, for example  by \cite{Artin 1971}, Proposition 1.4.
Clearly, the   $\Spec(R/\idealq)_\primid$ are the irreducible components of $\Spec(R_\primid)$.
The assertion now follows from the following lemma.
\qed

\begin{lemma}
\mylabel{henselian mimimal primes}
A local ring $A$ is henselian if and only if $A/\idealq$ is henselian for
every minimal prime ideal $\idealq\subset A$.
\end{lemma}

\proof
The condition is necessary by  \cite{EGA IVd}, Proposition 18.5.10.
For the converse, write $Y=\Spec(A)$. Let $f:X\ra Y$ be a finite morphism, and write $a_1,\ldots,a_n\in X$ for  the closed points. 
Since $X$ is quasicompact, each point $x\in X$ specializes to at least one closed point.
Since $Z=\overline{\left\{x\right\}}$  is irreducible, the composite morphism $Z\ra Y$ is finite, 
and $Y$ is henselian, it follows that $Z$ is local.
Consequently, we have a disjoint union $X=\Spec(\O_{X,x_1})\cup\ldots\cup\Spec(\O_{X,a_n})$.
Writing $X=\Spec(B)$, we conlcude that the maximal ideals $\maxid_i\subset B$ corresponding
to the closed points $a_i\in X$ are coprime, and thus $B=\prod B_{\maxid_i}$.
In turn, $X$ is a sum of local schemes.
\qed

\medskip
We now easily obtain the following useful criterion:

\begin{proposition}
\mylabel{criterion fppf-local}
A ring $R$ is fppf-local if and only if it is local henselian and
every finite fppf homomorphism $R\ra B$ admits a retraction.
\end{proposition}

\proof
According to Theorem \ref{amazing property}, the condition is necessary.
It is sufficient as well: Suppose that $R$ is local henselian, and
every finite fppf algebra admits a retraction. Let $R\ra C$
be an arbitrary fppf homomorphism. According to \cite{EGA IVd}, Corollary 17.16.2
there is a residue class ring $C/\ideala$ that is quasifinite and fppf over $R$.
Since $R$ is henselian, there is a larger ideal $\ideala\subset\idealb$ so that
$B=C/\idealb$ is finite and fppf. The latter admits, by assumptions,
a retraction $B\ra R$, and the composite map $C\ra B\ra R$ is the desired retraction of $C$.
\qed

\medskip
Now let $X$ be a scheme, and $P:(\Set) \ra X_\fppf$ be a point in the sense of topos-theory.
Applying the corresponding fiber functor $P^{-1}$ to the structure sheaf $\O_{X_\fppf}$,
we get a ring $\O_{X_\fppf,P}=P^{-1}(\O_{X_\fppf})$.

\begin{theorem}
\mylabel{stalks fppf-local}
Under the preceding assumptions, the ring $ \O_{X_\fppf,P}$ is fppf-local.
\end{theorem}

\proof
Choose a pro-object $(U_i)_{i\in I}$ in $(\fppf/X)$ so that the fiber functor is
of the form $P^{-1}(F)=\dirlim_{i\in I}\Gamma(U_i,F)$. Write $R_i=\Gamma(U_i,\O_{X_\fppf})$ 
and $R=\O_{X_\fppf,P}$, such that $R=\dirlim_{i\in I} R_i$.
According to Lemma \ref{affine pro-object} below, we may assume that the schemes $U_i$ are affine, in other words, $U_i=\Spec(R_i)$.

We first verify that the ring $R$ is local. 
In light of \cite{Hakim 1972}, Chapter III, Corollary 2.7, it suffices
to check that the ringed topos $(X_\fppf,\O_{X_\fppf})$ is locally ringed.
This indeed holds by loc.\ cit.\ Criterion 2.4, because for each $U\in(\fppf/X)$ and each $s\in\Gamma(U,\O_U)$,
the open subsets $U_s,U_{1-s}\subset U$ where $s$ respectively $1-s$ are invertible form a covering.

We next check that the local ring $R$ is henselian. Let $R\ra B$ be \'etale,
and suppose there is a retraction $B/\maxid\ra k$, where $\maxid\subset R$ is the maximal
ideal, and $k=R/\maxid$. We have to verify that this retraction extends to a retraction $B\ra R$,
compare \cite{EGA IVd}, Theorem 18.5.11.
Localizing $B$, we may assume that $B$ is local, such that we merely have to check that
there is a retraction $B\ra R$ at all.
According to \cite{EGA IVd} Proposition 17.7.8, there is an index $i\in I$ and some \'etale homomorphism $R_i\ra B_i$
with $B=B_i\otimes_{R_i}R$. Set $B_j=R_j\otimes_{R_i}B_i$ for $j\geq i$.
Invoking \cite{SGA 4a}, Expose IV, Section 6.8.7 again, we infer that
there is an index $j\geq i$ and some $R_i$-algebra homomorphism $R_j\ra B_i$, which gives
a retraction $B_i\otimes_{R_i}R_j\ra R_j$. This yields a direct system of retractions
$$
B_i\otimes_{R_i}R_k = B_i\otimes_{R_i}R_j\otimes_{R_j}R_k\lra R_{j}\otimes_{R_j}R_k=R_k.
$$
Passing to direct limits with respect to $k\geq j$ yields the desired retraction 
$B=\dirlim_{k\geq i} (B_i\otimes_{R_i} R_k)\ra R$.

We finally show that $R$ is fppf-local.
Let $R\ra B$ be a finite fppf homomorphism of rings. It suffices to check that it admits a retraction,
by Proposition \ref{criterion fppf-local}. Define $B_i$ as in the preceding paragraph.
There is an index $i\in I$ 
and an  homomorphism $R_i\ra B_i$ with $B=B_i\otimes_{R_i}R$, according to  \cite{EGA IVc}, Theorem 8.8.2.
Since the $R$-module $B$ is free of finite nonzero rank,
 we may assume that the same holds for the $R_i$-module $B_i$, by \cite{EGA IVc}, Corollary 8.5.2.5.
Set $V=U_i=\Spec(R_i)$, and $V'=\Spec(B_i)$. Again using \cite{SGA 4a}, Expose IV, Section 6.8.7
and arguing as above, one infers that the desired retraction $B\ra R$ exists.
\qed

\medskip
In the course of the preceding proof, we have used the following fact:

\begin{proposition}
\mylabel{affine pro-object}
Each topos-theoretical point $P:(\Set)\ra(\fppf/X)$ has a fiber functor isomorphic to  $F\mapsto\dirlim_{i\in I}\Gamma((U_i,p_i),F)$
for some pro-object $((U_i,p_i))_{i\in I}$ in $(\fppf/X)$ where all the $U_i$ are affine.
\end{proposition}

\proof
To simplify notation, write $\calC'=(\fppf/X)$  and consider the full subcategory $\calC\subset\calC'$ of objects
$(U,p)$ with $U$ affine, endowed with the induced Grothendieck topology.
Let $\shE'=X_\fppf$ and $\shE$ be the ensuing topoi of $\shU$-sheaves.
Given an object $(U,p)\in \calC'$, we denote by $U_\alpha\subset U$ the family of all 
affine open subschemes, and set $p_\alpha=p|U_\alpha$.
Clearly, $(U_\alpha,p_\alpha)\in\calC$ and $((U_\alpha,p_\alpha)\ra (U,p))_\alpha$ is a covering in $\calC'$. 
By the Comparison Lemma
(\cite{SGA 4a}, Expose III, Theorem 4.1), the restriction functor $\shE'\ra\shE$, $F\mapsto F|\calC$
is an equivalence of categories. In turn, every fiber functor on $\shE'$ is isomorphic to some
fiber functor coming  from a pro-object in the category $\calC$.
\qed

\section{Construction of fppf-local rings}
\mylabel{construction rings}

Let $R$ be a strictly local ring, that is, a henselian local ring with separably closed
residue field. Choose a universe $R\in\shU$  and some
ordinal $\sigma'\not\in\shU$.
Let $\sigma<\sigma'$ be the smallest ordinal that is not an element of  $\shU$. The goal of this
section is to attach, in a functorial way, 
a direct system $R_\nu\in\shU$ of strictly local rings, 
indexed by the well-ordered set
$$
\sigma=\left\{\nu\mid\text{$\nu$ ordinal with $\nu<\sigma$}\right\}
$$ 
of all smaller ordinals. The transition maps in this direct system will be faithfully flat and integral.
For each limit ordinal $\lambda<\sigma$, the local ring $R_\lambda$ will be fppf-local,
that is, every   fppf $R_\lambda$-algebra admits a retraction.
Maybe it goes without saying that all the rings $R_\nu$, $\nu>0$ are highly non-noetherian.

The construction of the rings is as follows:
Consider the   category $\calF=\calF(R)$ of finite fppf   $R$-algebras $A$ with $\Spec(A)$ connected.
Note that each such algebra is a fortiori isomorphic to some $R[T_0,\ldots,T_m]/(f_1,\ldots,f_r)$
for some integer $m\geq 0$ and some finite collection of polynomials $f_1,\ldots,f_r$.
Whence the set of isomorphism classes of objects in $\calF$
does not depend on the chosen universe, up to canonical bijection. Choose a  set 
$I=I(R)$ of such $R$-algebras, so that each isomorphism class is represented by precisely one
element of $I$.
Now consider the set $\Phi=\Phi(R)$ of all finite subsets of the set $I$,
endowed with the order relation coming from the inclusion relation. Clearly, the ordered set $\Phi$ is
filtered. Each of its elements $\varphi$ is  thus a finite set of certain finite fppf $R$-algebras.
Given an element $\varphi\in\Phi$, we form the finite fppf $R$-algebra
$$
A_\varphi=\bigotimes_{A\in\varphi}A.
$$
Here the tensor product denotes the unordered tensor product. 

Recall that for an   collection of $R$-modules
$(M_j)_{j\in J}$, indexed by some finite set $J$ of cardinality $n\geq 0$, the \emph{unordered tensor product} 
is the $R$-module of invariants
$$
\bigotimes_{j\in J}M_j=\left(\bigoplus_{\eta} M_{\eta(1)}\otimes\ldots\otimes M_{\eta(n)}\right)^{S_n}.
$$
Here the sum runs over all bijections $\eta:\left\{1,\ldots,n\right\}\ra J$, and the symmetric group $S_n$ acts from the
right on the sum
by permuting the summands:
$$
 (a_{\eta(1)}\otimes\ldots\otimes a_{\eta(n)})_\eta \cdot \sigma = 
(a_{\eta\sigma(1)}\otimes\ldots\otimes a_{\eta\sigma(n)})_{\eta\sigma}.
$$
Note that, for each choice of ordering $J=\left\{j_1,\ldots,j_n\right\}$,  
the obvious inclusion into the sum gives a canonical identification
$M_{j_1}\otimes\ldots\otimes M_{j_n}=\bigotimes_{j\in J}M_j$ of the ordinary tensor product with the unordered tensor product. 
However, the unordered tensor product has the advantage to be    functorial,  in the strict sense, with
respect to indexed $R$-modules $(M_j)_{j\in J}$. 
Here a morphism between $(M_j)_{j\in J}$ and $(M_k)_{k\in K}$ is given
by a map $m:J\ra K$ together with homomorphisms $M_j\ra N_{m(j)}$.

This functoriality ensures that $\varphi\mapsto A_\varphi$
is a direct system of $R$-algebras, and we already remarked that it is filtered. We denote by 
$$
R_+=\dirlim_{\varphi\in\Phi(R)} A_\varphi
$$ 
its direct limit. Clearly, all transition maps in this direct limit are finite fppf, so we may
regard each $A_\varphi$ as an $R$-subalgebra of $R_+$. 

\begin{lemma}
\mylabel{plus construction}
The ring $R_+$ is strictly local, and the homomorphism $R\ra R_+$ is flat, integral and local.
Moreover, we have $R_+\in\shU$, and $\Card(R_+)=\Card(R)$.
\end{lemma}

\proof
Using the notation from the beginning of this section, we start by checking 
that the   tensor products $A_\varphi\simeq A_1\otimes\ldots\otimes A_n$,
where $A_1,\ldots,A_n$ are the elements $\varphi$, are local. Obviously, 
$A_\varphi$ is a finite fppf $R$-algebra.
By definition, the schemes $\Spec(A_i)$ are connected. Whence $R\ra A_i$ are local maps of local rings, because $R$ is henselian. 
Furthermore, $A_i\otimes_Rk$, where
$k=R/\maxid_R$ is the residue field, is a finite local $k$-algebra. Their tensor product remains
local, because $k$ is separably closed. We conclude that there is a unique prime ideal in $A$
lying over $\maxid_R\subset R$. Since $\Spec(A)\ra\Spec(R)$ is a closed map, it follows that
$A$ is local, and that the map $R\ra A$ is local.
Passing to the filtered direct limit, the first assertion follows.

Since $R$, $A_\varphi$, $I(R)$ and whence $\Phi(R)$ are elements of the universe $\shU$,
the same must hold for the direct limit $R_+$.
By faithful flatness, the maps $A_\varphi\ra R_+$ are injective, whence $R_+=\bigcup_{\varphi\in \Phi} A_\varphi$ as a union of 
subrings that are finite fppf $R$-algebras. 
Now recall that $R$ is strictly local, and in particular infinite.
Let $\aleph_\iota$ be its cardinality.
It easily follows that each subring $A_\varphi\subset R_+$
and the index set $\Phi$ both have the same cardinality.  
Cardinal arithmetic  thus gives $\aleph_\iota\leq \Card(R_+)\leq \aleph_\iota\cdot\aleph_\iota =\aleph_\iota$.
Consequently $\Card(R_+)=\Card(R)$.
\qed

\medskip
The ring $R_+$, however, is never noetherian: 
If $n\geq 1$,  then $R[T]/(T^n)$ is  a finite fppf algebra with connected spectrum,  whence
isomorphic to some subring of $R_+$. It follows that there is an element $f\in R_+$ with $f^n=0$
but $f^{n-1}\neq 0$. In particular, the nilradical $\Nil(R_+)$ is not nilpotent.

Using \emph{transfinite recursion}, we now define a direct system of rings $R_\nu$, $\nu<\sigma$ as follows:
To start with,    set $R_0=R$. Suppose the direct system is already defined for all ordinals
 smaller that some $\nu<\sigma$. We then set
$$
R_{\nu}=\begin{cases}
(R_\gamma)_+ & \text{if $\nu=\gamma+1$ is a successor ordinal;}\\
\dirlim_{\gamma<\nu} R_\gamma & \text{if $\nu$ is a limit ordinal.}
\end{cases}
$$

\begin{proposition}
\mylabel{transfinite induction}
For each ordinal $\nu<\sigma$, the rings $R_\nu$ are strictly local, 
and the transition map $R_\gamma\ra R_{\nu}$, $\gamma\leq\nu$ are local,
faithfully flat, and integral.  
\end{proposition}

\proof
By transfinite induction. The assertion is trivial for $\nu=0$.
Now suppose that $\nu>0$, and that the assertion is true for all smaller ordinals.
If $\nu$ is a successor ordinal, the assertion follows from Lemma \ref{plus construction}.
If $\nu=\lambda$ is a limit ordinal, then $R_\lambda$ is a filtered direct limit
of strictly local ring with local transition maps, whence strictly local.
Moreover, the transition maps $R_\gamma\ra R_\lambda$ for $\gamma<\lambda$ are local, faithfully flat, and integral.
\qed

\begin{proposition}
\mylabel{dimension}
For each ordinal $\nu<\sigma$, we have $\dim(R_\nu)=\dim(R)$, and the residue field 
$k_\nu=R_\nu/\maxid_{R_\nu}$ is an algebraic closure of the residue field $k=R/\maxid_R$.
\end{proposition}

\proof
Since $R\subset R_\nu$ is integral and faithfully flat, the first statement follows
from \cite{AC 8-9}, Chapter VIII, \S 2, No.\ 3, Theorem 1. 
As to the second assertion, the field extension $k\subset k_\nu$ is clearly algebraic.
Let $P\in R[T]$ be a monic polynomial,
consider the finite fppf $R$-algebra $A=R[T]/(P)$, and let $A_\maxid$ be the localization
at some maximal ideal $\maxid\subset A$. Then $P$ has a root in $A$, and 
 $A_\maxid$ is isomorphic to 
a subring of $R_\nu$. Whence we have a homomorphism $A\ra k_\nu$. 
It follows that each monic polynomials with coefficients in $k$
has a root in $k_\nu$.
\qed

\begin{proposition}
\mylabel{element of universe}
For each ordinal $\nu<\sigma$, the ring  $R_\nu$ is   an element of the chosen universe $\shU$.
\end{proposition}

\proof
By transfinite induction. The assertion is obvious for $\nu=0$.
Now suppose that $\nu>0$, and that the assertion holds for all smaller ordinals.
If $\nu=\gamma+1$ is a successor ordinal, then $R_\nu=(R_\gamma)_+\in\shU$ by Proposition \ref{plus construction}.
If $\nu$ is a limit ordinal, then $R_\nu=\dirlim_{\gamma<\nu} R_\gamma\in\shU$ because the 
$R_\gamma$ and the index set, which equals the set $\nu$, are members of the universe $\shU$.
\qed

\begin{theorem}
\mylabel{fppf acyclic}
For each limit ordinal $\lambda<\sigma$, the ring $R_\lambda$ is  fppf-local. 
\end{theorem}

\proof
Given a fppf homomorphism $R\ra B$, we have to show that it admits a retraction.
It suffices to treat the case that $B$ is finite fppf, according to Proposition \ref{criterion fppf-local}.
Since $B$ is of finite presentation, there is some ordinal $\nu<\lambda$ and
some $R_\nu$-algebra $B_\nu$ 
with $B\simeq B_\nu\otimes_{R_\nu}R_\lambda$ (see \cite{EGA IVb}, Lemma 5.13.7.1). Moreover, we may assume that 
$B_\nu$ is finite and fppf. Since tensor products commute with filtered direct limits,
the canonical map $\dirlim_\gamma(B_\nu\otimes_{R_\nu}R_\gamma)\ra B$ is bijective,
where the direct limit runs over all ordinals $\nu\leq\gamma<\lambda$.

Choose $A_\nu\in I(R_\nu)$ and an isomorphism of $R_\nu$-algebras $h:B_\nu\ra A_\nu$.
Consider the singleton $\varphi=\left\{A_\nu\right\}\in\Phi(R_\nu)$, such that $A_\varphi=A_\nu$
in the notation introduced above. By the very definition of $R_{\nu+1}$, there exists an $R_\nu$-algebra
homomorphism $B_\nu\stackrel{h}{\ra} A_\varphi\ra R_{\nu+1}$, which gives a retraction $B_\nu\otimes_{R_\nu}R_{\nu+1}\ra R_{\nu+1}$.
Tensoring with $R_\gamma$ over $R_{\gamma+1}$, $\gamma\geq \nu+1$ we get a direct system of retractions
$$
B_\nu\otimes_{R_\nu}R_\gamma=B_\nu\otimes_{R_\nu}R_{\nu+1}\otimes_{R_{\nu+1}}R_\gamma\lra 
R_{\nu+1}\otimes_{R_{\nu+1}}R_\gamma=R_\gamma.
$$
Passing to direct limits yields the desired retraction $B\ra R_\lambda$.
\qed

\medskip
For later use, we record the following fact:

\begin{lemma}
\mylabel{finite subalgebras}
For each ordinal $\nu<\sigma$ and each finite subset $S\subset R_\nu$, there
is an $R$-subalgebra $B\subset R_\nu$ containing $S$ so that the
structure map $R\ra B$ is finite and fppf.
\end{lemma}

\proof
By transfinite induction. The case $\nu=0$ is trivial.
Now suppose that $\nu>0$, and that the assertion is true for all smaller ordinals.
If $\nu$ is a limit ordinal, then $R_\nu=\dirlim_{\gamma<\lambda}R_\gamma$,   
there is some ordinal $\gamma<\nu$
with $S\subset R_\gamma$, and the induction hypothesis, together with flatness
of $R_\gamma\subset R_\nu$, yields the assertion.

Now suppose that $\nu=\gamma+1$ is a successor ordinal.
Write $R_\nu=\dirlim A_\varphi$ as a filtered union of finite   fppf local $R_\gamma$-subalgebras 
and choose some index $\varphi $ so that $S\subset A_\varphi$.
Since $R_\gamma$ is local, the underlying $R_\gamma$-module of $A_\varphi$ is free.
The same holds for $A_\varphi/R_\gamma$, because the unit element $1\in A_\varphi$
does not vanish anywhere. Thus we may extend $b_1=1$ to an $R_\gamma$-basis $b_1,\ldots, b_m\in A_\varphi$,
and write
$$
s=\sum_k c_{sk}b_k\quadand b_i\cdot b_j=\sum_kc_{ijk}b_k
$$
for some coefficients $c_{sk},c_{ijk}\in R_\gamma$, where $s\in S$ and $1\leq i,j,k\leq m$.
Form the finite subset $S'=\left\{c_{sk},c_{ijk}\right\}\subset R_\gamma$ comprising all these coefficients.
By induction hypothesis, there is a finite fppf $R$-subalgebra $B'\subset R_\gamma$
containing $S'$. Now consider the canonical $R$-linear map
$$
\bigoplus_{i=1}^mB'b_i\lra R_{\nu}.
$$
This map factors over $A_\varphi\subset R_{\nu}$, and it is injective, because the $b_1,\ldots,b_m\in A_\varphi$ are 
$R_\gamma$-linearly independent.  Let $B\subset R_{\nu}$ be its image.
By construction, $S\subset B$, and $B'\subset B$ is a direct summand of free $B'$-modules of finite rank,
in particular, an fppf ring extension. It follows that $R\subset B$ is finite fppf.
\qed

\medskip
Given $\nu<\sigma$, consider the set of all $R$-subalgebras $B_i\subset R_\nu$, $i\in I_\nu $ so 
that the structure map $R\ra B_i$ is finite fppf. We regard $I_\nu$ is an ordered set,
where the order relation is the inclusion relation.
 
\begin{proposition}
\mylabel{index filtered}
The ordered set $I_\nu$ is filtered, each $B_i$ is a local $R$-algebra such that the structure
map $R\ra B_i$ is local, finite and fppf, and $R_\nu=\bigcup_{i\in I_\nu} B_i$.
\end{proposition}

\proof
If follows from Lemma \ref{finite subalgebras} that $R_\nu$ is the union of the $B_i$.
To see that the union is filtered, let $B_i,B_j\subset R_\nu$ be two such
subrings.  Let $S_i\subset B_i$ be an $R$-basis, and similarly $S_j\subset B_j$.
Then $S=S_i\cup S_j$ is a finite subset of $R_\nu$, and Lemma \ref{finite subalgebras} gives
us the desired subalgebra $B\subset R_\nu$ containing $B_i$ and $B_j$.

Each $B_i$ is by definition finite fppf over $R$. Since $R$ is strictly local,
it remains to check that $\Spec(B_i)$ is connected.
Since $B_i\ra R_\nu$ is injective, the image of the continuous map
$\Spec(R_\nu)\ra\Spec(B_i)$ contains every generic point.
Obviously, $B_i\subset R_\nu$ is integral, whence the continuous map is surjective.
Since $\Spec(R_\nu)$ is connected, so must be its continuous image $\Spec(B_i)$.
\qed

\medskip
The direct system $R_\nu$, $\nu<\sigma$ is functorial: Suppose that $f:R\ra R'$ is a local homomorphism
between strictly local rings. With the notation introduced at the beginning of this section,
we get a functor
$$
\calF(R)\lra\calF(R'),\quad A\longmapsto A'=A\otimes_RR',
$$
and thus induced maps of ordered sets $I(R)\ra I(R')$ and $\Phi(A)\ra\Phi(A')$, $\phi\mapsto\phi'$.
The latter are not necessarily injective, but in any case induce morphisms
$A_\phi\ra A_{\phi'}$. In turn, we get a natural homomorphism of direct limits
$R_+\ra R'_+$. Using transfinite induction, one finally obtains the desired homomorphism
$R_\nu\ra R'_\nu$ of direct systems. One easily checks that this is functorial.

\section{Points in the fppf topos}
\mylabel{construction points}

Let $X$ be a scheme. Choose a universe $X\in\shU$, and let $\sigma$ be the smallest ordinal
not contained in this universe.
Given a geometric point $\bar{a}:\Spec(\Omega)\ra X$ and a limit ordinal $\lambda<\sigma$,
we call
$$
\O_{X,\bar{a},\lambda}=(\O_{X,\bar{a}})_\lambda
$$
the \emph{fppf-local ring attached to the geometric point and the limit ordinal}, as defined in Section \ref{construction rings}.
The goal now is to construct  a point  $P=P_{\bar{a},\lambda}:(\Set)\ra X_{\fppf}$
in the sense of topos-theory, together with a canonical identification 
$$
\O_{X_\fppf,P}=P^{-1}(\O_{X_\fppf}) = \O_{X,\bar{a},\lambda}.
$$
Actually, the isomorphism class of $P_{\bar{a},\lambda}\in\Points(X_\fppf)$ depends only on the image point $a\in X$ of the geometric point $\bar{a}$,
and will be denoted by $P_{a,\lambda}\in |X_\fppf|$.
This will give  a continuous section $a\mapsto P_{a,\lambda}$ for the canonical map
$|X_\fppf|\lra |X_\Zar|=X$
of topological spaces.

The main step is the construction of a pro-object in $(\fppf/X)$ that will be  isomorphic
to the pro-object of neighborhood for the topos-theoretical point $P_{\bar{a},\lambda}$.
Our first task to to find a suitable index category for such a pro-object:

Fix a geometric point $\bar{a}:\Spec(\Omega)\ra X$ and some ordinal $\nu<\sigma$. For the moment,
this can be either a limit ordinal or a successor ordinal.
We now define the index category $I_{\bar{a},\nu}$ as follows: The objects
are 5-tuples
$$
(V_0,V_1,\phi,U,\psi),
$$
where $a\in V_0\subset X$ is an affine open neighborhood, $V_1$ is an affine \'etale $V_0$-scheme, 
$\phi:\Spec(\O_{X,\bar{a}})\ra V_1$
is a morphism, $U$ is a finite fppf $V_1$-scheme,
and $\psi:\Spec(\O_{X,\bar{a},\nu})\ra U$ is a  morphism.
We demand that the diagram 
\begin{equation}
\label{neighborhood diagram}
\begin{CD}
@.                    \Spec(\O_{X,a}) @<\can<< \Spec(\O_{X,\bar{a}}) @<\can<< \Spec(\O_{X,\bar{a},\nu})\\\
@.                    @V\can VV                 @V\phi VV                  @VV\psi V\\
X                @<<< V_0             @<<< V_1                   @<<< U
\end{CD}
\end{equation}
is commutative, and that the  resulting morphism of affine schemes
\begin{equation}
\label{schematically dominant}
\Spec(\O_{X,\bar{a},\nu}) \lra U\times_{V_1}\Spec(\O_{X,\bar{a}})
\end{equation}
is \emph{schematically dominant}, that is, induces an injection on global sections of the structure sheaf.
The morphisms
$$
(V_0',V_1',\phi',U',\psi')\lra (V_0,V_1,\phi,U,\psi)
$$
in the category $I_{\bar{a},\nu}$ are   3-tuples $(h_0,h_1,h)$, 
where $h_0:V_0'\ra V_0$ and
$h_1:V_1'\ra V_1$ and $h:U'\ra U$ are morphisms of schemes.
We demand  that the diagram 
\begin{equation}
\label{big commutative}
\begin{gathered}
\xymatrix{
 &\Spec(\O_{X,a})\ar[d]\ar[ddr]  	& \Spec(\O_{X,\bar{a}})\ar[d]_{\phi'}\ar[l]\ar[ddr]^/-2em/{\phi} & \Spec(\O_{X,\bar{a},\lambda})\ar[d]_{\psi'}\ar[l]\ar[ddr]^/-2em/\psi\\
X\ar[dr]_{\id}& V_0'\ar[dr]_{h_0}\ar[l]		& V_1'\ar[l]\ar[dr]_{h_1}			& U'\ar[l]\ar[dr]_h\\
 & X				& V_0\ar[l]						& V_1\ar[l]	& U\ar[l]
}
\end{gathered}
\end{equation}
is commutative. Note that $h_0$ is an inclusion between the two open subschemes $V_0',V_0\subset X$,
and $h_1$ is a refinement   between the \'etale neighborhoods $V_1',V_1\ra X$ of the geometric point $\bar{a}$.

Given an object $(V_0,V_1,\phi,U,\psi)$, the composite morphism $U\ra V_1\ra V_0\subset X$
is quasifinite and fppf, whence we may regard $U$ as an object in $(\fppf/X)$.
This yields a covariant functor $I_{\bar{a},\nu}\ra(\fppf/X)$, which on morphism is defined
as $(h_0,h_1,h)\mapsto h$. 
The corresponding contravariant functor $I_{\bar{a},\nu}^\op\ra(\fppf/X)$
is actually a pro-object, which means the following: 

\begin{proposition}
\mylabel{category filtered}
The   opposed category $I_{\bar{a},\nu}^\op$ is filtered.
\end{proposition}

\proof
Working with $I_{\bar{a},\nu}$ rather then the opposed category, we have to check
two things: First, for any two given objects, there is some object and morphisms from it to the given objects.
Second,   any two morphisms with the same domain and range become equal after composing with some morphism
from the right.

We start with the former condition: Suppose  
\begin{equation}
\label{two objects}
(V_0',V_1',\phi',U',\psi') \quadand  (V_0'',V_1'',\phi'',U'',\psi'').
\end{equation}
are two objects. Choose an affine open neighborhood $a\in V_0\subset X$
contained in $V_0'\cap V_0''$. Base-changing the data to $V_0$ over $X$, we easily
reduce to the case that $V_0'=V_0=V_0''$.
Similarly, we may assume $V_1'=V_1''$ and $\phi'=\phi''$.

To simplify notation, now write 
$$
R=\O_{X,\bar{a}} \quadand R_\nu=(\O_{X,\bar{a}} )_\nu 
$$
Let $B'\subset R_\nu$ be the images of the ring of global sections
with respect to the morphism of affine schemes $\Spec(R_\nu)\ra U\times_{V_1}\Spec(R)$.
The definition of objects in $I_{\bar{a},\nu}$ ensures that $\Spec(B')=U\times_{V_1}\Spec(R)$.
Define $B''\subset R_\nu$ analogously. Then $B',B''$ are finite fppf as $R$-algebras,
and $R$-subalgebras inside $R_\nu$.
According to Proposition \ref{index filtered}, they are contained in some larger finite fppf $R$-subalgebra $B\subset R_\nu$.
Since $R=\O_{X,\bar{a}}$ can be regarded as  the filtered direct limit of the global section rings of the \'etale neighborhoods
of the geometric point
$$
\Spec(\Omega)\subset\Spec(R)\stackrel{\phi'=\phi''}{\lra} V_1'=V_1'',
$$ we find some \'etale neighborhood $V_1\ra V_1'=V_1''$  
so
that $\Spec(B)\ra\Spec(R)$ arises via base-change from some finite fppf scheme $U\ra V_1$.
Passing to smaller \'etale neighborhoods, we may assume that the inclusion maps $B',B''\subset B$ inside $R_\nu$ 
of finite fppf $R$-algebras are
induced by some $V_1$-morphisms $U\ra U'$ and $U\ra U''$.
Let $\phi:\Spec(R)\ra V_1$ and $\psi:\Spec(R_\nu)\ra U$ be the canonical morphisms.
Then $(V_0,V_1,\phi,U,\psi)$ is an object in $I_{\bar{a},\nu}$ and by construction  
has morphisms to both of the given objects in (\ref{two objects}).

It remains to verify the second condition. Suppose we have two arrows
$$
\xymatrix{
(V_0',V_1',\phi',U',\psi') \ar@<2pt>[r] \ar@<-2pt>[r] & (V_0,V_1,\phi,U,\psi)
}
$$
called   $(h_0,h_1,h)$ and $(k_0,k_1,k)$. We have to show that they become equal after
composing from the right  with some morphism. According to the commutative diagram (\ref{big commutative}),
both  $h,k:U'\ra U$ are $V_1$-morphisms. Consider the two morphisms  
$$
U'\times_{V_1}\Spec(\O_{X,\bar{a}})\stackrel{h,k}{\lra} U\times_{V_1}\Spec(\O_{X,\bar{a}})\stackrel{\pr}{\lra} U.
$$
coming form base-change. These    coincide, because  the morphism in (\ref{schematically dominant}) is
schematically dominant.
Since $U',U$ are  of finite presentation over $V_1$, 
the morphisms 
$$
U'\times_{V_1}V_1''\stackrel{h,k}{\lra} U\times_{V_1}V_1''\stackrel{\pr}{\lra} U
$$
coming from  base-change of $h$ and $k$ to some \'etale neighborhood $V_1''\ra V_1$ of  the geometric point
$\phi:\Spec(\Omega)\ra V_1$ become
identical (\cite{EGA IVc}, Theorem 8.8.2). Choosing the neighborhood small enough, we may
assume that it factors over $V_1'$. Let $\phi'':\Spec(\O_{X,\bar{a}})\ra V_1''$ be the canonical map.
Define $U''=U'\times_{V_1}V_1''$ and let $\psi'':\Spec(\O_{X,\bar{a},\lambda})\ra U''$
be the canonical map. 
The resulting morphism
$(V_0',V_1'',\phi'',U'',\psi'')\ra (V_0',V_1',\phi',U',\psi')$ does the job.
\qed

\medskip
We now  have a pro-object
$$
I_{\bar{a},\nu}\lra(\fppf/X),\quad (V_0,V_1,\phi,U,\psi)\longmapsto U
$$
and obtain a covariant functor
\begin{equation}
\label{inverse image functor}
X_\fppf\lra (\Set),\quad F\longmapsto \dirlim\Gamma(U,F),
\end{equation}
where the direct limit runs over all objects $(V_0,V_1,\phi,U,\psi)\in I_{\bar{a},\nu}$ 
This functor respects finite inverse limits, because the opposite of the index category is filtered.
In the special case $F=\O_{X_\fppf}$, the morphisms $\psi:\Spec(\O_{X,\bar{a},\nu})\ra U$ induce a canonical 
homomorphism
\begin{equation}
\label{map in direct limit}
\dirlim\Gamma(U,\O_{X_\fppf})\lra \O_{X,\bar{a},\nu}
\end{equation}
of rings.

\begin{proposition}
\mylabel{stalk correct}
The preceding homomorphism of rings (\ref{map in direct limit}) is bijective.
\end{proposition}

\proof
The map in question is surjective: Set $R=\O_{X,\bar{a}}$ and $R_\nu=(\O_{X,\bar{a}})_\nu$.
According to Proposition \ref{index filtered}, every element $c\in R_\nu$ is contained in some 
$R$-subalgebra $C\subset R_\nu$ so that the homomorphism $R\ra C$ is finite fppf.
Write $R=\dirlim_{i\in I} R_i$ as a filtered direct limit with \'etale neighborhoods  $\Spec(R_i)\ra X$
of the geometric point $\bar{a}:\Spec(\Omega)\ra X$.
For some index $j\in I$, there is a finite fppf $R_j$-algebra $C_j$ with $C_j\otimes_{R_j}R=C$.
Then 
$$
C=C_j\otimes_{R_j} R=C_j\otimes_{R_j}\dirlim_{i\geq j} R_i=\dirlim_{i\geq j} (C_j\otimes_{R_j}R_i).
$$
Whence there is some index $i\geq j$ so that $c\in C$ lies in the image  $C_i=C\otimes_{R_j}R_i$.
Replacing $j$ by $i$ and $C$ by $C_i$, we thus may assume that $c\in C$ is in the image of of $C_j$.

Set $U=\Spec(C_j)$ and $V_1=\Spec(R_j)$, 
and let $\phi:\Spec(R)\ra V_1$ and $\psi:\Spec(R_\nu)\ra U$ be the canonical morphisms.
The image of the structure map $V_1\ra X$, which is \'etale, is an open neighborhood of $a\in X$,
whence contains some affine open neighborhood $a\in V_0$. Base-changing 
with $V_0$, we may assume that $V_1\ra X$ factors over $V_0$. Replacing $V_1$ by an affine
open neighborhood of the image of $\phi$, we may again assume that $V_1$ is affine.
The upshot is that the tuple $(V_0,V_1,\phi,U,\psi)$ is an object of the index category $I_{\bar{a},\nu}$.
By construction, the element $c\in R_\nu$ lies in the image of $\Gamma(U,\O_{X_\fppf})$.

The map  is injective as well: Suppose we have an object $(V_0,V_1,\phi,U,\psi)\in I_{\bar{a},\nu}$
and some local section $s\in\Gamma(U,\O_{X_\fppf})=\Gamma(U,\O_U)$ whose image in $R_\nu$ vanishes.
By definition of the index category, the map
$$
\Gamma(U,\O_U)\otimes_{\Gamma(V_1,\O_{V_1})} R=\Gamma(U\times_{V_1}\Spec R,\O_{U\times_{V_1}\Spec R})\lra R_\nu
$$
is injective, whence $s\otimes 1$ vanishes in the left hand side.
Now regard $R=\dirlim R_i$ as a filtered direct limit for \'etale neighborhoods
$\Spec(R_i)\ra V_1$ of the induced geometric point $\phi:\Spec(\Omega)\ra V_1$. Then it vanishes
already in $\Gamma(U,\O_U)\otimes_{\Gamma(V_1,\O_{V_1})} R_i$ for some index $i\in I$.
Set $V_1'=\Spec(R_i)$ and $U'=U\otimes_{\Gamma(V_1,\O_{V_1})} R_i$, endowed
with the canonical morphisms $\phi':\Spec(R)\ra V_1'$ and $\psi':\Spec(R_\nu)\ra U'$.
Then $(V_0,V_1',\phi',U',\psi)$ is an object in $I_{\bar{a},\nu}$, endowed with
a canonical morphism to $(V_0,V_1,\phi,U,\psi) $, on which the pullback
of $s$ vanishes.
\qed

\begin{theorem}
\mylabel{topos-theoretical point}
Suppose our $\nu=\lambda$ is a limit ordinal.
Then there is a topos-theoretical point $P_{\bar{a},\lambda}:(\text{\rm Set})\ra X_\fppf$
whose inverse image functor $P^{-1}_{\bar{a},\lambda}$ equals the functor in (\ref{inverse image functor}).
\end{theorem}

\proof
We apply the criterion given in \cite{SGA 4a}, Expose IV, 6.8.7.
Let $W\in(\fppf/X)$ be an object,   $(W_\alpha\ra W)_{\alpha\in\Lambda}$   an fppf covering, 
$(V_0,V_1,\phi,U,\psi)\in I_{\bar{a},\lambda}$   an index, and $U\ra W$
be a  $X$-morphism. We have to find a larger index $(V'_0,V'_1,\phi',U',\psi')$,
some $\alpha\in \Lambda$ and  a  morphism $U'\ra W_\alpha$ making the diagram
$$
\begin{CD}
U' 	@>>> 	W_\alpha\\
@VVV 		@VVV\\
U	@>>>	W
\end{CD}
$$
commutative. 
To simplify notation, set $R=\O_{X,\bar{a}}$ and $R_\lambda=\O_{X,\bar{a},\lambda}$.
Consider the base-changes $\Spec(R_\lambda)\times_WW_\alpha$. Choose some $\alpha\in \Lambda$
so that the closed point in the local scheme $\Spec(R_\lambda)$ is in the image of the projection.
By flatness, the projection is thus surjective.
According to Theorem \ref{stalks fppf-local}, there is a morphism
$\Spec(\O_{X,\bar{a},\lambda})\ra W_\alpha$ making the diagram
$$
\xymatrix{
						&		& W_\alpha\ar[d]\\
\Spec(\O_{X,\bar{a},\lambda})\ar[r]\ar[rru]	& U\ar[r]	& W
}
$$
commute. Using Proposition \ref{stalk correct}, together with \cite{EGA IVc}, Theorem 8.8.2,
we conclude that there is some morphism $(V'_0,V'_1,\phi',U',\psi')\ra (V_0,V_1,\phi,U,\psi)$ 
and a $W$-morphism $U'\ra W_\alpha$, as desired.
\qed

\medskip
Let $a\in X$ be the image of the geometric point $\bar{a}$.
If $\bar{b}$ is another geometric point on $X$ whose image points $b$ equals $a$,
there is a $\kappa(a)$-isomorphism $\kappa(\bar{a})\ra \kappa(\bar{b})$,
which comes from a unique isomorphism of strictly local rings $\O_{X,\bar{a}}\ra\O_{X,\bar{b}}$.
By functoriality, it extends to an isomorphism $\O_{X,\bar{a},\nu}\ra\O_{X,\bar{b},\nu}$,
which finally  yields an isomorphism of inverse systems
$$
(V_0,V_1,\phi, U, \psi)\longmapsto (V_0,V_1,\phi f, U, \psi  f).
$$
We conclude that the isomorphism class   of the topos-theoretical point $P_{\bar{a},\lambda}\in\Points(X_\fppf)$
only depends on the image point $a\in X$, and we write this isomorphism class as $P_{a,\lambda}\in |X_\fppf|$. 

\begin{theorem}
\mylabel{continuous section}
Let $\lambda<\sigma$ be a limit ordinal. Then the map 
$$
X=|X_\Zar|\lra |X_\fppf|,\quad a\longmapsto P_{a,\lambda}
$$
is  section for the canonical projection $|X_\fppf|\ra |X_\Zar|=X$.
\end{theorem}

\proof
Suppose that $V\subset X$ is an open subscheme
that is a neighborhood of the topos-theoretical point $ P_{\bar{a},\lambda}$.
Then there is an index $(V_0,V_1,\psi,U,\phi)\in I_{\bar{a},\lambda}$
having an $X$-morphism $U\ra V$. By the commutative diagram (\ref{neighborhood diagram}),
the diagram
$$
\begin{CD}
\Spec(\O_{X,\bar{a},\lambda})	@>>> \Spec(\O_{X,a})\\
@V\phi VV 			     @VVV\\
U 				@>>> X
\end{CD}
$$
is commutative as well. In turn, we have $a\in V$. According to \cite{SGA 4a}, Expose IV, Section 7.1
there is a unique point $a'\in X$ so that the open subschemes of $X$ that are
neighborhoods of $P_{\bar{a},\lambda}$ are   neighborhoods of $a'$. It follows that $a'=a$.
Hence $a\mapsto P_{\bar{a},\lambda}$ is a section.
\qed

\medskip
In light of Proposition \ref{fppf homeomorphisms} and Lemma \ref{bijective locales}, the existence
of a section now yields our main result:

\begin{theorem}
\mylabel{identifiaction sober}
The continuous map $|X_\fppf|\ra |X_\Zar|$ induces an identification 
$|X_\fppf|_\sob=|X_\Zar|_\sob=X$ of sober spaces.
\end{theorem}


\end{document}